\documentstyle[twoside,amssymb,amsmath,12pt]{article}
\def\C{{\Bbb C}}

\def\R{{\Bbb R}}
\def\P{{\Bbb P}}
\def\Z{{\Bbb Z}}
\def\A{{\Bbb A}}
\def\G{{\Bbb G}}
\def\T{{\Bbb T}}

\newtheorem{prop}{Proposition}[subsection]
\newtheorem{dfn}[prop]{Definition}
\newtheorem{theo}[prop]{Theorem}

\newtheorem{rem}[prop]{Remark}
\newtheorem{coro}[prop]{Corollary}
\newtheorem{lem}[prop]{Lemma}

\newtheorem{abs}[prop]{$\,$}

\newcommand{\npo}{{\P}^1}
\newcommand{\npt}{{\P}^2}
\newcommand{\npr}{\P}
\newcommand{\npop}{{\P}^1 \times {\P}^1}
\newcommand{\FP}{\rm Fano polyhedron $P$}

\title{\sc On the Classification of Toric Fano $4$-folds} 

\author{{\sc Victor V. Batyrev} \\
\small  {\em Mathematisches Institut, Universit\"at T\"ubingen}   \\
\small  {\em Auf der Morgenstelle 10,  72076  T\"ubingen, Germany}  \\
\small  {\em e-mail: batyrev@bastau.mathematik.uni-tuebingen.de} \\
 }

\begin{document}

\date{}

\maketitle

\begin{abstract}
The biregular classification of smooth $d$-dimensional 
toric Fano varieties of dimension $d$ is equivalent to the classification 
of special simplicial polyhedra $P$  in $\R^d$, so called 
Fano polyhedra, up to an isomorphism of the standard 
lattice $\Z^d\subset \R^d$. 
In this paper we explain the complete biregular classification of 
all $4$-dimensional smooth toric Fano varieties. The main result 
states that there  exist exactly 123 different types of toric Fano 
$4$-folds up to isomorphism. 
\end{abstract}

\tableofcontents

\newpage

\section{Introduction}

A smooth projective $d$-dimensional algebraic variety $V$ over $\C$ 
is called a 
Fano manifold if the anticanonical sheaf of $V$ is ample. In the case  
$d=2$ the surface  
$V$ is also called a  Del Pezzo surface. It is a classical result that 
all Del Pezzo  surfaces can be obtained from $\P^2$ and 
$\P^1 \times \P^1$  
by blowing ups of respectively $r \leq 8$ and $r \leq 7$  
points in general position 
(see e.g. \cite{manin}). This shows that there exist exactly $10$ different 
types of Del Pezzo surfaces up to deformation. 

First fundamental  results towards the complete classification of Fano 
$3$-folds up to deformation were received by Iskovskih \cite{isk1,isk2} 
more than 20 years ago. The complete classification of Fano $3$-folds has 
been obtained later in the papers of Mori and Mukai 
\cite{MM1,MM2,MM3} (see also \cite{murre,cutkosky}).  
It has been proved by Koll\'ar, Miyaoka and Mori \cite{KMM1,KMM2} 
that there exists only finitely many smooth Fano manifolds 
of fixed dimension $d$ 
up to deformation. For  Fano $d$-folds with Picard number one this 
statement was independently obtained by Nadel \cite{nadel}). 

The complete classification of all Fano $4$-folds up to deformation still 
remains an open problem. One expects that the complete list of Fano $4$-folds 
up to deformation could include  many hundreds of  different types. 
For this reason, it looks more meaningful to restrict the classification 
of Fano manifolds of dimension $d \geq 4$ to some special classes for 
which the complete list is no so large.

In the  present paper we consider Fano $d$-folds which are 
simultaniously toric varieties, i.e., admitting a regular effective 
action of a $d$-dimensional algebraic torus $(\C^*)^d$
\cite{dan2,fulton,oda}. Since toric varieties are determined by some discrete 
combinatorial data,  they do not admit nontrivial algebraic deformations. 
For this reason, it make sense to classify toric Fano $d$-folds up to 
{\em biregular isomorphism}. Toric Fano manifolds 
of dimension $d \leq 3$ have been classified  in the author's paper 
\cite{bat1} and independently in the paper of  
K.$\,$Watanabe$\,$\&$\,$M.$\,$Watanabe 
\cite{wat}. There exist exactly $5$ different toric Del Pezzo surfaces and 
exactly $18$ different toric Fano $3$-folds. 
Some special classes of $d$-dimensional Fano manifolds of  
arbitrary dimension $d$  have been classified by 
Voskresenski\^i$\,$\&$\,$Klyachko \cite{vosk.klyachko} and 
Ewald \cite{ewald}.  

It is important to stress that the classification problem of 
toric Fano $d$-folds up to isomorphism can be reformulated purely into a  
combinatorial classification problem of special convex 
polyhedra $P \subset \R^d$, so 
called Fano polyhedra, up to linear unimodular transformation from 
$GL(d, \Z)$. Every $d$-dimensional Fano polyhedron $P$ determines a toric 
Fano $d$-fold $V(P)$. Fano polyhedra form a  special subclass of  
{\em reflexive 
polyhedra} which were introduced in \cite{bat.dual}.  
We want to remark  that  an explicit   classification  of  higher 
dimensional Fano polyhedra $P$ up to unimodular transformations needs 
a special method   for  describing such  polyhedra. For instance, if  
we described a $d$-dimensional Fano polyhedron $P$ with $n$ vertices 
just by the $d\times n$-matrix $M(P)$ consisting of the   
coordinates of the  vertices, then  
we would  meet the following two difficulties: 

i) the $d\times n$-matrix $M(P)$ doesn't show  
the combinatorial structure of $P$ which contain a lot of information 
about the geometry of the corresponding toric Fano $d$-fold $V(P)$; 

ii) having  two $d\times n$-matrices  
$M(P_1)$ and  $M(P_2)$ which describe given Fano polyhedra $P_1$  
and  $P_2$,  it is not easy  to decide whether $P_1$ and $P_2$ 
are equivalent up to an unimodular linear transformation or not. 

Both these difficulties disappear if one uses the language of 
{\em primitive collections} and {\em primitive relations} for 
describing toric Fano $d$-folds. This language was introduced in 
\cite{bat5} and developed in \cite{bat.cox} in the connection with 
\cite{cox}.

The main result of the paper is the complete classification of 
toric Fano $4$-folds up to isomorphism. This classification 
was obtained in the author's PhD thesis \cite{bat4}. 
In Diplomarbeit \cite{evertz} S. Evertz independently 
has classified $4$-dimensional Fano having  $8$ vertices using 
some computer calculations. Both papers \cite{bat4} and \cite{evertz} 
have many common ideas and both use the combinatorial classification 
of simplicial $4$-dimensional polyhedra with $8$ vertices 
due to  Grunbaum and Streedharan \cite{grun}.  Comparing the list of  
$4$-dimensional Fano polyhedra with $8$ vertices in  \cite{bat4} and 
\cite{evertz}, some discrepancy in results has been observed: only $3$ 
polyhedra of the combinatorial type $P^8_{26}$ 
(instead of $5$  \cite{bat4})  were  found in \cite{evertz}, and
$2$  Fano polyhedra of the combinatorial type 
$P^8_{28}$ from \cite{evertz} were missing in the list of  
\cite{bat4} (the combinatorial type $P^8_{28}$ was 
excluded from considerations in \cite{bat4} by some error).

The classification of toric Fano $4$-folds was used by Y. Nakagawa 
in \cite{nakagawa1,nakagawa2} 
for finding all toric Fano $4$-folds which admit  
Einstein-K\"ahler metrics. This result of Nakagawa extends the one of 
Mabuchi \cite{mabuchi} for toric Fano $3$-folds.   
\bigskip

\noindent
{\bf Acknowledgement.} I would like to thank Professors
G\"unter Ewald, Tadao Oda and Toshiki  
Mabuchi for their interest to the results in  
\cite{bat4} and stimulating conversations. I am also grateful to 
Yasuhiro Nakagawa who pointed out me on some misprints \cite{nakagawa0} in 
an earlier version of the present  paper.

\section{Toric Fano manifolds of dimension $d$}

\subsection{Fano polyhedra} 

\begin{dfn} {\rm Let $P$ be a convex polyhedron in ${\R}^d$,   
${\cal V}(P) = \{ v_1, \ldots, v_n \}$ the set of all vertices 
of $P$. We call  $P$ a  {\bf  Fano polyhedron}   
if  the  following conditions are satisfied: 

 (i) the elements of ${\cal V}(P)$ belong to the 
standard integral lattice ${\Z}^d  
\subset  {\R}^d$;

(ii) $P$ contains the lattice element $(0, \ldots, 0) \in {\Z}^d$  in its 
interior;

(iii) $P$ is a simplicial polyhedron, i.e.,   each  face  of  $P$  
is  a  simplex;

(iv) vertices $v_{i_1} , \ldots , v_{i_d}$ of 
any  $(d-1)$-dimensional 
face $F = [v_{i_1} , \ldots , v_{i_d}] $ 
of $P$ form a ${\Z}$-basis of the lattice 
${\Z}^d$, i.e., 
the coordinates of $v_{i_1} , \ldots , v_{i_d}$  form a matrix $A$ 
with ${\rm det} \, A = \pm 1$.}
\label{def.fano}
\end{dfn}

\begin{dfn} 
{\rm Let $P$ be a Fano polyhedron  and 
$\{ v_{i_1} , \ldots , v_{i_m} \}$  be the set of 
vertices of a $(k-1)$-dimensional 
face of $F \subset P$ $(1 \leq m \leq d)$. 
We denote by  $\sigma (F)$ the $m$-dimensional 
cone consisting of all nonnegative ${\R}$-linear combinations of 
$v_{i_1} , \ldots , v_{i_m}$, i.e., 
\[ \sigma (F) := \{ \lambda_1 v_{i_1} + \ldots + \lambda_m v_{i_m} \in \R^n \;
: \; \lambda_i \geq 0\; 
(1 \leq i \leq m) \} . \] 
The system of  cones 
\[ \Sigma(P) = \{ 0 \} \cup \{ \sigma(F) \}_{F \subsetneq P}, \]
where $F$ runs over all proper faces of $P$, we call the {\bf polyhedral fan}
associated  with the Fano polyhedron $P$. }
\end{dfn}

\begin{dfn}{\rm  Let $P$ be a Fano polyhedron. A 
subset 
\[ {\cal P} = \{ v_{i_1} , \ldots , v_{i_k} \} \subset {\cal V}(P) \]  
will be  called  a {\bf  primitive  collection}  if  it 
satisfies the conditions

(i) ${\cal P}$ does not belong to a  cone from $\Sigma (P)$;

(ii) any proper subset of ${\cal P}$ is contained in some 
cone from $\Sigma (P)$.}
\label{prim.coll}
\end{dfn}

\begin{dfn}
{\rm Let $P$ be a Fano polyhedron and ${\cal P}=  \{ v_{i_1} , \ldots , 
v_{i_k} \} $ 
be a primitive collection of its vertices. 
Denote by $\sigma ({\cal P})$ 
the  cone of { minimal dimension} in $\Sigma(P)$  
containing the  integral 
point $s({\cal P}) = v_{i_1} + \ldots + v_{i_k}$ 
(Since the cones of the fan $\Sigma (P)$ cover the whole space ${\R}^d$, 
the point $s({\cal P})$ must belong to at least one  cone 
$\sigma(F) \subset \Sigma (P)$ for some face $F \subset P$.)
Let $v_{j_1} , \ldots , v_{j_m}$ are generators  
of $\sigma({\cal P})$. Then the  element $s({\cal P})$ has a unique 
representation as a 
positive integral linear combination of $v_{j_1} , \ldots , v_{j_m}$ :
\[ s({\cal P}) =  c_1 v_{j_1 } + \cdots + c_k v_{j_m}, \; \; c_i > 0, \; \; 
c_i \in {\Z}. \]
(If $s({\cal P}) =0$, then the set $\{v_{j_1} , \ldots , v_{j_m} \}$ is 
assumed to be empty.)  
We will call the linear  relation 
\[ R({\cal P}) :\;  
 v_{i_1} + \ldots + v_{i_k} - c_1 v_{j_1}  - \ldots - c_m y_{j_m} = 0 \] 
among the vertices of $P$ the {\bf  primitive 
relation} corresponding to the primitive collection ${\cal P}$. The 
positive integers $c_1, \ldots, c_m$ we  call the {\bf coefficients} of the 
primitive relation  $R({\cal P})$. The integer 
\[ \Delta({\cal P}) := k - \sum_{i=1}^{m} c_i \]
we call the {\bf  degree} of the primitive collection ${\cal P}$.}
\end{dfn}

Our next purpose is to show  that primitive relations are  useful 
for combinatorial descriptions  of Fano polyhedra $P$.  

\begin{dfn} 
{\rm Let $P_1$  and $P_2$ be two Fano  polyhedra.  
Then $P_1$ and  $P_2$ are called {\bf  combinatorically equivalent} 
if and only if there exists a bijective 
 mapping 
$$\rho\; : \; {\cal V}(P_1) \rightarrow {\cal V}(P_2)$$  
which  
respects the face-relation, i.e. , $v_{i_1} , \ldots , v_{i_k}$ are vertices 
of a $(k-1)$-dimensional face of $P_1$ if and only if 
$\rho(v_{i_1}) , \ldots , \rho(v_{i_k})$ are vertices of 
a $(k-1)$-dimensional face of $P_2$.}
\end{dfn}

By Definition \ref{prim.coll}, one immediatelly 
obtains the following characterization of the combinatorial equivalence 
class of $P$ using  primitive collections:

\begin{prop}  Let $P_1$  and $P_2$ be two Fano  polyhedra.  
Then $P_1$ and  $P_2$ are   combinatorically equivalent 
if and only if there exists a bijective 
 mapping 
$$\rho\; : \; {\cal V}(P_1) \rightarrow {\cal V}(P_2)$$ 
 which
 induces a one-to-one correspondence between primitive collections 
of vertices in ${\cal V}(P_1)$ and ${\cal V}(P_2)$.  
\label{comb-type}
\end{prop}

\begin{dfn} {\rm  Two Fano polyhedra $P_1$ and $P_2$ are called 
{\bf  isomorphic} if there exists an authomorphism $\rho \in 
{\rm GL}(d, {\Z})$ of the lattice ${\Z}^d$ such that 
$\rho (P_1) = P_2$. }
\end{dfn}
 
Assume that two Fano polyhedra $P_1$ and $P_2$  are isomorphic. Then the 
induced   bijection $\rho\, :\, {\cal V}(P_1) \rightarrow {\cal V}(P_2)$ 
between vertices obviously respects linear relations among them. 
It is less evident to see that the opposite  statement holds, so that we 
have the following:   

\begin{prop}
Two Fano polyhedra $P_1$ and  $P_2$ are isomorphic if and only if 
there exist a bijective 
mapping  between ${\cal V}(P_1)$ and ${\cal V}(P_2)$  which respects  
not only primitive collections of vertices in ${\cal V}(P_1)$ and 
${\cal V}(P_2)$, but also the corresponding primitive 
relations. 
\label{isomorphism} 
\end{prop} 

The less evident part "if"   
in \ref{isomorphism} follows immediately  from the following 
statement which we formulate without proof: 

\begin{lem}
Let $P$ be a $d$-dimensional Fano polyhedron having $n$ vertices. 
Denote by $L(P)$ sublattice of rank $n-d$ in ${\Z}^n$ consisting 
of all linear relations with integral coefficients among elements 
of ${\cal V}(P)$. 
Then $L(P)$ is generated by primitive relations. 
\label{lemma1}
\end{lem}

\begin{prop} Let  $P$ be  a Fano polyhedron and ${\cal P} 
= \{ v_{i_1} , \ldots , v_{i_k} \} \subset {\cal V}(P)$ 
is a  primitive collection. Then $\Delta({\cal P}) >0$. 
In particular,   exist only  finitely many 
possibilities for the coefficients $c_i$ in the  primitive relation 
$R({\cal P})$ if $P$ has a fixed dimension $d$.   
\label{coeff-prim}
\end{prop}

\noindent
{\em Proof. }
It follows from \ref{prim.coll} that  
the rational point 
\[ r({\cal P})= \frac{1}{k} \left(  v_{i_1} + \ldots + v_{i_k}\right)  \]
belongs to the interior of $P$. Therefore,   
the corresponding primitive  relation  $R({\cal P})$ has 
positive   integer   coefficients  $\{ c_{i} \}$   
satisfying the inequality 
\[ c_1 + \ldots + c_m < k.\] 
Obviously, $k \leq {\rm dim} \; P + 1 \; = \; d+1$. This 
implies that there exists only a finite number of possibilities 
for positive integers 
$c_i$ if the dimension $d$ is fixed. 
\hfill $\Box$

\begin{prop}
The number of vertices $n = n(P)$ of any  $d$-dimensional Fano polyhedron 
$P$  is not greater than $2(2^d -1)$. 
\label{vertices}
\end{prop}

\noindent
{\em Proof. }
Consider the canonical surjective homomorphism 
\[ \alpha\; : \; {\Z}^d \rightarrow {\Z}^d/(2{\Z})^d. \]
First we  remark  that no  vertex $v$ of $P$ belongs 
to the kernel of $\alpha$. Otherwise all coordinates of $v$ would be 
divisible by $2$. This contradicts to \ref{def.fano}(iv). 

Assume now that $\alpha(v_i) = \alpha(v_j) \neq 0$ for  
two different vertices $v_i, v_j \in P$. Then $(v_i + v_j)/2$ 
is an integral point belonging to $P$. Since the only lattice points belonging 
to the boundary of $P$ are its vertices and 
$(v_i + v_j)/2$ can not be a vertex of $P$ (see  \ref{def.fano}(iv)), 
we conclude that $(v_i + v_j)/2 =0$. 
Thus $v_i$ and $v_j$ are centrally symmetric vertices 
of $P$. This implies that preimage $\alpha^{-1}(x)$ of any  nonzero 
element $x \in {\Z}^d/ (2{\Z})^d$ contains at most two 
different vertices of $P$. Since ${\Z}^d/ (2{\Z})^d$ contains $2^d -1$ nonzero 
elements, 
we obtain the required  upper bound  for the number $n(P)$. 
\hfill $\Box$ 

\begin{rem} {\rm The above exponential upper bound  for the number of 
vertices of $P$ is sharp only for $d=1,2$. A better polynomial upper bound  
\[ n(P) \leq d^2 +1, \; {\rm for}\; d > 2 \]
has been obtained by 
Klyachko and Voskresenski\^i in \cite{vosk.klyachko}.}
\end{rem}

\begin{theo} Let $d$ be a fixed positive integer. Then there  
exist  only  finitely   many  $d$-dimensional Fano polyhedra  
up to isomorphism.
\label{finiteness}
\end{theo}

\noindent
\noindent {\em Proof.} By \ref{comb-type} and \ref{vertices}, there exist 
only finitely many different 
combinatorial types of $d$-dimensional Fano polyhedra. On the other hand, 
for a fixed combinatorial type of $P$ and for a fixed  primitive 
collection ${\cal P} \subset {\cal V}(P)$ there exist 
only finitely many possibilities for the corresponding primitive 
relation $R({\cal P})$ (see \ref{coeff-prim}).  Therefore there exist only 
finitely many possibilities for  primitive relations among vertices 
of a Fano polyhedron of a fixed dimension $d$. 
It remains to apply \ref{isomorphism}. \hfill $\Box$

\begin{rem} {\rm Theorem \ref{finiteness} was proved in 
\cite{bat2,borisovs,vosk.klyachko} by other methods.}
\end{rem}

\subsection{Toric Fano $d$-folds} 

Let $P$ be a $d$-dimensional Fano polyhedron.  
It follows from the theory of toric varieties \cite{dan2,fulton,oda} 
that the fan $\Sigma(P)$ defines a smooth projective toric variety $V(P)$ 
having ample anticanonical sheaf, i.e., a toric Fano $d$-fold. For 
convenience, we give an explicit geometric description of the toric 
Fano $d$-fold  $V(P)$ using 
our language of primitive collections and 
primitive relations:

\begin{dfn}
{\rm 
Let ${\A}^n$ be  $n$-dimensional affine space over  
${\C}$ with the complex coordinates 
$z_1, \ldots, z_n$, where $n$ is the number of vertices of $P$. We establish 
the one-to-one correspondence $z_i \leftrightarrow v_i$ between 
the coordinates $z_1, \ldots, z_n$ and the vertices $v_1, \ldots, v_n$. 
If  ${\cal P} = \{ v_{i_1}, \ldots, v_{i_k} \} 
\subset {\cal V}(P)$ be a primitive collection, then  we define 
an   affine subspace ${\A}({\cal P}) \subset {\A}^n$: 
\[ {\A}({\cal P}) := \{ (z_1, \ldots, z_n) \in {\A}^n\;: \; 
z_{i_1} = \cdots = z_{i_k} = 0 \}. \]
Denote by $U(P)$ the Zariski open subset 
\[ U(P):= {\A}^n \setminus \bigcup_{{\cal P} \subset {\cal V}(P)} 
{\A}({\cal P}). \]} 
\end{dfn} 

\begin{dfn} 
{\rm Let  $L(P) \subset {\Z}^n$ be the 
$(n-d)$-dimensional sublattice consisting of 
integral linear relations among vertices of $P$, i.e., 
\[ L(P): = \{ (\lambda_1, \ldots, \lambda_n) \in {\Z}^n\; : \; 
\lambda_1 v_1 + \ldots + \lambda_n v_n = 0 \}. \] 
For every element $\lambda = (\lambda_1, \ldots, \lambda_n)$ we define 
a $1$-parameter subgroup ${T}_{\lambda}$ in $({\C}^*)^n$ acting 
canonically on 
${\A}^n$ as follows 
\[ t (z_1, \ldots, z_n) = (t^{\lambda_1}z_1, \ldots, t^{\lambda_n}z_n), \; 
t \in {\C}^*. \] 
Denote by ${T(P)}$ the $(n-d)$-dimensional algebraic 
subgroup in $({\C}^*)^n$ generated by all $1$-parameter subgroups 
${T}_{\lambda}$, $\lambda \in L(P)$. } 
\end{dfn} 

\begin{dfn}
{\rm We define the toric manifold $V(P)$ as the quotient of $U(P)$ modulo 
the canonical linear action of ${T(P)}$ on $U(P)$ (it is easy to see 
that $T(P)$ acts free on $U(P)$). }
\end{dfn}

\begin{theo} The toric manifold 
$V(P)$ is a smooth projective  $d$-dimensional toric variety 
with ample anticanonical sheaf.  Any smooth compact $d$-dimensional  
toric variety $V$  having ample anticanonical sheaf, i.e., 
a toric Fano manifold,  is isomorphic to 
$V(P)$ for some $d$-dimensional Fano polyhedron $P$. 
Moreover, two  Fano  varieties  $V(P_1)$  
and $V(P_2 )$ corresponding to two $d$-dimensional Fano 
polyhedra $P_1 $  and $P_2 $  are  biregular 
isomorphic (as abstract algebraic varieties) 
if  and  only  if  $P_1$   and  $P_2$  are  isomorphic (as Fano polyhedra).
\end{theo}

\noindent
{\em Proof}. We explain only the way how an isomorphism of 
$V(P_1)$ and $V(P_2)$ as 
abstract complex manifolds implies isomorphism of 
the corresponding polyhedra 
$P_1$ and $P_2$. The rest part of the proof is an easy 
consequence of the standard theory of toric varieties 
\cite{dan2,fulton,oda}. 

Let $\psi \;: \; V(P_1) \rightarrow  V(P_2)$ be a biregular 
isomorphism of 
two toric Fano $d$-folds. Let $\G_i:= Aut(V(P_i))$ $(i=1,2)$ be the 
algebraic group of biregular authomorphisms of $V(P_i)$. 
Denote by $\T_i \subset \G_i$ $( i=1,2)$ the   
maximal $d$-dimensional torus canonically embedded as 
open subset in $V(P_i)$. The isomorphism $\psi$ induces an 
isomorphism 
$$\tilde{\psi}\; : \; \G_1 \to \G_2 $$ 
\[ \gamma \mapsto \psi \circ \gamma \circ \psi^{-1}, \;\; 
\gamma \in \G_1. \]
Now we show that there exists an isomorphism $\varphi \,: \, 
V(P_1) \rightarrow  V(P_2)$ such that the corresponding isomorphism 
$\tilde{\varphi}\, : \, \G_1 \to \G_2$ has the property 
$\tilde{\varphi}({\T}_1) = {\T}_2$. Indeed,  
$\tilde{\psi}({\T}_1)$ and ${\T}_2$  are two maximal 
tori in the linear algebraic group $\G_2$. 
By the well-known theorem of Borel \cite{borel}, 
${\T}_2$ and $\tilde{\psi}({\T}_1)$ are conjugate 
by some element $\gamma_0 \in \G_2$, i. e. $\tilde{\psi}({\T}_1) = 
\gamma^{-1}_0 \T_2 \gamma_0$.  Now we define $\varphi := 
\gamma_0 \circ \psi$. By our definition, we have  
$$\tilde{\varphi}(\T_1) = \gamma_0 (\tilde{\psi}(\T_1) \gamma_0^{-1} = 
\T_2.$$  Therefore 
$\varphi$ is an isomorphism of toric varieties $V(P_1)$ and 
$V(P_2)$ which respects the torus actions, i.e., an equivariant isomorphism. 
It follows from the theory of toric varieties that 
any equivariant isomorphism 
of toric varieties is determined by an isomorphism of fans $\Sigma(P_1)$ 
and $\Sigma(P_2)$. This isomorphism determines  an isomorphism of 
Fano polyhedra $P_1$ and $P_2$. 
\hfill $\Box$

\subsection{Primitive relations and extremal rays of Mori}  

The theory of extremal rays introduced by Mori \cite{mori1,mori2} 
was the main 
technical tool in the classification of Fano $3$-folds \cite{MM1,MM2,MM3}. 
It is natural to 
expect a similar role of extremal rays in the classification of 
toric Fano $d$-folds. The language of primitive relations 
turns out to be very convenient for this purpose:  

\begin{theo} 
The group $A_1(V(P))$ of $1$-dimensional cycles on $V(P)$ relative to 
numerical equivalence is canonically isomorphic to the $(n-d)$-dimensional 
lattice $L(P)$ of integral relations among vertices of $P$. 
Moreover,   the cone of Mori $\overline{NE}(V(P))  \subset 
A_1(V(P))\otimes {\R}$ is canonically isomorphic to the cone 
generated by all primitive 
relations  
\[ v_{i_1} + \cdots + v_{i_k} - c_1 v_{j_1} - 
\cdots - c_m v_{j_m} = 0 \]
which we consider as elements of 
$L(P) \subset {\Z}^n$. 
\label{mori.cone}
\end{theo}

\begin{rem} 
{\rm By theorem \ref{mori.cone},  every $1$-dimensional 
face of $\overline{NE}(V(P))$, e.g., an extremal ray,  
is generated by a primitive relation. However, the converse  is not true 
in general. There might be primitive relations which corresponds 
interior lattice points on  faces  of dimension $\geq 2$ in  
$\overline{NE}(V(P))$.}
\end{rem}

\begin{theo}
Assume that a primitive relation 
\[ R({\cal P}) \;: \;v_{i_1} + \cdots + v_{i_k} - 
c_1 v_{j_1} - \cdots - c_m v_{j_m} = 0 \] 
defines a generator $\lambda_{\cal P}$ of some $1$-dimensional face of 
$\overline{NE}(V(P))$ (e.g., this condition is  always satisfied 
if $\Delta({\cal P})=1)$. 
Then the following statements hold:

{\rm (i)} The vertices $v_{j_1},  \ldots , v_{j_m}$ together with 
any $k - 1$ vertices from $\{ v_{i_1}, \ldots, v_{i_k} \}$ are 
generators of a $(m + k -1)$-dimensional cone in $\Sigma(P)$. 

{\rm (ii)} A $1$-dimensional cycle on $V(P)$ representing
the class of $R({\cal P})$ in $A_1(V(P))$ is a smooth projective rational  
curve with the normal bundle 
\[ \underbrace{{\cal O}(1) \oplus \cdots \oplus {\cal O}(1)}_{k-2} \oplus 
\underbrace{{\cal O} \oplus \cdots \oplus {\cal O}}_{d+1-k-m}
 \oplus {\cal O}(-c_1) \oplus 
\cdots \oplus {\cal O}(-c_m). \]
In particular, the anticanonical degree of $\lambda_{\cal P} \in A_1(V(P))$ 
equals $\Delta({\cal P})$. 

{\rm (iii)} If $m=1$ $($i.e., $v_{i_1} + \cdots + v_{i_k} - 
cv_{j}$ is a primitive relation$)$, then the exceptional locus 
of the extremal contraction corresponding to $\lambda_{\cal P}$ is a 
divisor $D_j$ on 
$V(P)$ which  is a 
${\P}^{k-1}$-bundle over a smooth $(d -k)$-dimensional toric 
variety $W$. 
\label{ex-ray2}
\end{theo}

\begin{coro} Let $P$ be a Fano $d$-polyhedron, and let $V$ 
be the corresponding toric Fano $d$-fold. Assume that there exists 
a primitive relation $R({\cal P})$ among vertices of $P$ as follows 
\[ v_i  + v_j - v_l =0 . \]
Then $R({\cal P})$ gives rise to an extremal ray with the  extremal 
contraction $V \rightarrow V'$ to  a smooth projective toric variety $V'$ 
$($$V'$ is not necessary a Fano $d$-fold in general$)$. 
\label{except}
\end{coro} 

\begin{rem} 
{\rm The statements of \ref{mori.cone} and \ref{ex-ray2} can be obtained 
from the Mori's thery for projective toric varieties (see \cite{reid} or 
\cite{oda}, 2.5).} 
\end{rem}   

The language of primitive collections and primitive relations 
among generators of $1$-dimensional cones in $\Sigma$  can be used 
in description of arbitrary (not necessary Fano) smooth compact 
toric varieties ${\P}_{\Sigma}$ (see \cite{bat5}): 

\begin{prop}
Let ${\P}_{\Sigma}$ be a compact smooth $d$-dimensional toric variety 
corresponding to a complete $d$-dimensional regular fan $\Sigma$ with 
generators $\{ v_1, \ldots, v_n \}$. Then  
the anticanonical class of ${\P}_{\Sigma}$ is ample $($resp. numerically 
effective$)$ if and 
only if for every primitive relation 
\[  v_{i_1} + \cdots + v_{i_k} - 
c_1 v_{j_1} - \cdots - c_m v_{j_m} = 0 \]
one has $k - \sum_{i=1}^{m} > 0$ $($resp. $k - \sum_{i=1}^{m} \geq 0$ $)$. 
\label{crit}
\end{prop}

There exist another  combinatorial method of description of smooth compact 
$d$-dimensional  toric manifolds via so called {\bf  weighted triangulations} 
of $(d-1)$-dimensional sphere $S^{d-1}$ \cite{oda}. 
This method gives an explicit 
information about all $1$-dimensional toric strata and their normal bundles. 
For toric Fano manifolds the  degrees of normal bundles to all 
$1$-dimensional strata must be at least $ -1$. This  requirement 
puts  some restrictions of the number of faces of Fano polyhedra $P$: 

\begin{theo} 
Let $f_i(P)$ be the number of $i$-dimensional faces of a $d$-dimensional 
Fano polyhedron $P$. Then 
\[ 12f_{d-3}(P) \geq (3d-4) f_{d-2}(P).\] 
Moreover, the equality holds if and only if all $1$-dimensional torus 
invariant strata on $V(P)$ are rational curves 
having  the anticanonical degree $1$, i.e., the degree of the normal 
bundle to every $1$-dimensional stratum on $V(P)$ is $-1$. 
\label{weights-1}
\end{theo}

\noindent 
{\em Proof.} Let $w(P)$ be the sum of all weights in the weighted 
triangulation of $S^{d-1}$ defining the fan $\Sigma(P)$. 
It is know that a sequence $w_1, \ldots, w_s$ of integers  
in a weighted circular graph defining a $2$-dimensional smooth  
projective variety satisfies the condition (see \cite{oda}. p.45):
\[ \sum_{j=1}^s a_j = 12 - 3s. \]
Using this formula for each $2$-dimensional toric stratum in $V(P)$ 
corresponding to a $(d-3)$-dimensional face of $P$, we obtain   
\[ w(P)= 12 f_{d-3}(P) - 3(d-1)f_{d-2}(P). \]
On the other hand, all $1$-dimensional toric strata in $V(P)$ are 
parametrized by $(d-2)$-dimensional faces of $P$. Since the sum of 
weights by every $(d-2)$-dimensional face of $P$ equals the degree
of the normal bundle, we conclude
\[  w(P) \geq  -f_{d-2}(P).\] 
Moreover, the last inequality becomes equality iff the degree of the normal 
bundle to every $1$-dimensional stratum on $V(P)$ is $-1$. This implies the 
statement of theorem. 

\hfill $\Box$

\begin{rem} 
{\rm In \cite{evertz} the number $w(P)= 
12 f_{d-3}(P) - 3(d-1)f_{d-2}(P)$ was 
called the {\bf total weight} of $P$. } 
\end{rem}  

\subsection{Projections of Fano polyhedra}

\begin{dfn}
{\rm Let $P$ be a $d$-dimensional Fano polyhedron, $v_i \in {\cal V}(P)$ 
a vertex of $P$, ${\R}\langle v_i \rangle$ the $1$-dimensional subspace 
in ${\R}^n$ generated by $v_i$. We denote by 
$$\pi_i\; : \; {\R}^n \rightarrow 
{\R}^n/{\R}\langle v_i \rangle$$ 
the canonical epimorphism. 
The image $P_i = \pi_i(P)$ 
we call a {\bf  $\pi_i$-projection} of $P$. }
\end{dfn}

\begin{rem} 
{\rm It would be perfect for the classification of $d$-dimensional 
Fano polyhedra by induction on $d$, if the polyhedra 
$P_i$ were again  Fano polyhedra of dimension $d-1$. 
Unfortunately, this is not true in general. However,  the polyhedron 
$P_i$ is always very close to a Fano polyhedron.} 
\end{rem} 
 
\begin{prop}
The polyhedron $P_i$ is the convex hull of all points $\pi_i(v_j)$ such 
that the segment $[ v_i, v_j]$ is an edge  of $P$. Moreover, $0 \in 
{\R}^n/{\R}\langle v_i \rangle$ is the unique interior lattice point of 
$P_i$.
\label{conv-h} 
\end{prop}

\noindent 
{\em Proof.} By definition, $P_i$ is the convex hull of $\pi_i$-images 
of vertices of $P$. Let $P_i' \subset P_i$ be the the convex hull of 
$\pi_i$-images of all vertices $v_j$ such that $[v_i,v_j]$ is a face 
of $P$. It remains to prove that $P_i \subset P_i'$. 

Assume that $[ v_i, v_k ]$ is not an edge  of $P$. Then 
$\{ v_i, v_k \}$ is a primitive collection. There are two possibilities 
for the corresponding primitive relation:

{\sl Case} I. $v_i + v_k = 0$. This implies that $\pi_i(v_k) =0$. 
Obviously, $0 \in P_i'$. 

{\sl Case} II. $v_i + v_k = v_l$. By \ref{ex-ray2}(iii), 
$[ v_i, v_l ]$ is a face of 
$P$. So $\pi_i(v_k) = \pi_i(v_l) \in P_i'$. 

It follows from \ref{def.fano}(iii) that $P_i$ is covered 
by $\pi_i$-images of all $(d-1)$-faces of $P$ containing $v_i$. Hence, by 
\ref{def.fano}(iv), $0$ is the unique interior lattice point of $P_i$.
\hfill $\Box$

\begin{prop}
The anticanonical class  of toric divisor $D_i \subset V(P)$ corresponding 
to the vertex $v_i$  is always numerically effective. In particular, 
$P_i$ is always a reflexive polyedron. 
\label{num-eff} 
\end{prop}

\noindent 
{\em Proof.} Let $\alpha_i\,: \, 
{\rm Pic}\,V(P) \rightarrow {\rm Pic}\, D_i$ be 
the natural surjective mapping induced by restriction of Cartier 
divisors. This induces the injective mapping of dual lattices  
$\alpha^*_i\,: \, A_1(D_i)  
\rightarrow A_1(V(P))$. Since the cone of Mori of any toric variety 
is generated by classes $1$-strata \cite{reid}, 
$\alpha^*_i(\overline{NE}(D_i))$ is a face 
of the cone $\overline{NE}(V(P))$.
 So it is sufficient to prove  that for  any primitive 
relation $R({\cal P})$ such that $\lambda_{\cal P} 
\in \alpha^*_i(\overline{NE}(D_i))$. Thus,  for a primitive relation 
\[ R({\cal P})\;: \;   
v_{i_1} + \ldots + v_{i_k} - 
c_1 v_{j_1} - \ldots - c_m v_{j_m} = 0 \] 
representing a $1$-statum in $D_i$ (see \ref{mori.cone}), one has 
$\langle  -K_{D_i}, \lambda_{\cal P} \rangle \geq 0$. 
By the adjunction formula, we have  
$$\langle  -K_{D_i}, \lambda_{\cal P} \rangle = 
\langle  -K_{V(P)},  \lambda_{\cal P} \rangle - \langle D_i,  
\lambda_{\cal P} \rangle.$$   
Since $-K_{V(P)}$ is ample, $\langle  -K_{V(P)},  \lambda_{\cal P} 
\rangle \geq 1$. 
By \ref{ex-ray2}(ii), we obtain
\[ \langle D_i,  \lambda_{\cal P}  \rangle = \left\{  \begin{array}{l} 
1\;\; {\rm if}\; i \in \{ i_1, \ldots, i_k \}, \\
-c_l\;\; {\rm if}\; i= j_l, \\
0\;\; {\rm otherwise,} \end{array} \right. \]
i.e. ,  $\langle D_i,  \lambda_{\cal P}  \rangle \geq -1$. 
Hence, $\langle  -K_{D_i},  \lambda_{\cal P} 
\rangle \geq 0$. \hfill $\Box$
  
\begin{coro}
Assume that all primitive collection containing $v_i$ have degree $\Delta 
\geq 2$. 
Then $D_i$ is a $(d-1)$-dimensional toric Fano manifold, i.e., $P_i$ is 
a $(d-1)$-dimensional Fano polyhedron. 
\end{coro}

\begin{coro}
Let $\Sigma_i(P) \subset {\R}^n/{\R}\langle v_i \rangle$ 
be the $(d-1)$-dimensional complete regular fan defining $D_i$, $\Gamma$ 
a $(d-2)$-dimensional face of $P_i$. Then there exist only the following 
possibilities for $\Gamma$: 

{\rm (i)} $\Gamma$ contains $(d-1)$ lattice points which generate a 
$(d-1)$-dimensional cone from $\Sigma_i(P)$; 

{\rm (ii)} There exist a   primitive collection ${\cal P} = \{v,  v_{i_1}, 
\ldots, v_{i_k} \}$ of degree $1$, i.e, a primitive relation 
\[  v_i + v_{i_1} + \cdots + v_{i_k} - 
c_1 v_{j_1} - \cdots - c_m v_{j_m} = 0, \;\; k = \sum_{s=1}^{m} c_s, \]
such that $\Gamma$ contains 
exactly $d$ lattice points satisfying the uniquely determined 
linear relation    
\[ \pi_i(v_{i_1}) +  \cdots + \pi_i(v_{i_k}) -  c_1\pi_i(v_{j_1}) - 
\cdots -  c_m \pi_i(v_{j_m})=0. \]
\label{cor-gamma}
\end{coro} 

\begin{dfn}
{\rm A nonzero lattice point  $v \in P_i$ is called {\bf double point} 
if there exist 
two different vertices $v_j, v_k \in P$ such that $\pi_i(v_j) = 
\pi_i(v_k) = v$. }
\end{dfn}

\begin{prop}
Let $\Gamma$ be a $(d-2)$-dimensional face of $P_i$. 

{\rm (i)} If  $\Gamma$ contains $d$ lattice points, then no one of these 
point is double. 

{\rm (ii)} If $\Gamma$ contains $(d-1)$ lattice points, then $\Gamma$ 
contains at most $1$ double point. 
\label{commun}
\end{prop}

\noindent 
{\em Proof.} It is clear that $P$ has a supporting  
$(d-1)$-dimensional affine hyperplane $H \subset {\R}^d$ such 
that $\pi_i(H)$ is a $(d-2)$-dimensional affine hyperplane 
in ${\R}^n/{\R}\langle v_i \rangle$ with  
$\pi_i(H) \cap P_i = \Gamma$. Now both statements (i) and (ii) 
follow from the fact that $H$ can  not contain more than $d$ vertices 
of $P$. 
\hfill $\Box$. 

\begin{prop}
Assume that $P$ contains two centrally symmetric vertices $v_i$ and $v_j$, 
i.e., $v_i + v_j = 0$ is a primitive relation. If 
\[ R({\cal P}) \;:\; v_i + v_{i_1} + \cdots + v_{i_k} - 
c_1 v_{j_1} - \cdots - c_m v_{j_m} = 0  \]
is a primitive relation and  $k = \sum_{s=1}^m c_s$,  then 
$c_1 = \cdots = c_m =1, \; k =m$ 
and 
\[ v_j + v_{j_1} + \cdots + v_{j_m} - 
 v_{i_1} - \ldots -  v_{i_k} = 0  \]
is again a primitive relation. 
\label{cent-sym1}
\end{prop} 

\noindent 
{\em Proof.} By \ref{ex-ray2}(i), the set of 
vertices 
$\{  v_{i_1},  \ldots,  v_{i_k}, v_{j_1}, \ldots, v_{j_m} \}$ 
is a set of vertices  of a $(k+m)$-dimensional face $F$ of $P$. In particular, 
one has $k + m \leq d$. Since the lattice points 
$\pi_i(v_{i_1}), \ldots, \pi_i(v_{i_k}), \pi_i(v_{j_1}),  
\ldots, \pi_i(v_{j_m})$ satisfy the linear relation  
 \[ \pi_i(v_{i_1}) +  \cdots + \pi_i(v_{i_k}) -  c_1\pi_i(v_{j_1}) - 
\cdots -  c_m \pi_i(v_{j_m}) =0,  \]
there exists a $(d-2)$-dimensional face $\Gamma$ 
of $P_i$ containing $\pi_i(F)$. 
Since $v_i$ and $v_j$ are centrally symmetric, we can identify two 
projections $\pi_i(P)$ and $\pi_j(P)$. Applying \ref{cor-gamma} for 
the $(d-2)$-dimensional face $\Gamma$ of $P_j = \pi_j(P)$, we conclude
that there must a primitive relation
\[ R({\cal P}')\; : \; 
v_j + v_{l_1} + \cdots + v_{l_r} - d_1 v_{m_1} - \cdots - 
d_s v_{m_s} = 0, \]
where $\pi_j(v_{l_1}), \ldots, \pi_j(v_{l_s}), \pi_j(v_{m_1}), \ldots, 
\pi_j(v_{m_s})$ are lattice points in $\Gamma$ satisfying 
the linear relation  
\[ \pi_j(v_{l_1}) +  \cdots + \pi_j(v_{l_r}) -  d_1\pi_j(v_{m_1}) - 
\cdots -  d_s \pi_j(v_{m_s}) =0.  \] 
By \ref{cor-gamma}, there exists a unique linear relation among 
all $d$ lattice points of $\Gamma$. This implies that either
$k=r$, $m=s$ and $c_i =d_i$ $(i = 1, \ldots, m)$, or 
$k=s$, $m=r$ and $c_i=1$ $(i = 1, \ldots, m)$. The first case 
is impossible, because the sum of $ R({\cal P})$ and
$ R({\cal P}')$ would be a 
nontrivial linear relation among vertices of $F$.  
The second case implies  exactly the required statement. 
\hfill $\Box$

\begin{prop} 
Let $v \in P_i$ be a nonzero 
double point, i.e., $0 \neq v = \pi_i(v_j) = \pi_i(v_k)$ 
for some two different vertices $v_j, v_k \in P$. 
Then exactly one of the following
two  situations holds: 

(i) $v_i + v_k = v_j$ is a primitive relation; 

(ii) $v_i + v_j = v_k$ is a primitive relation. 
\label{double-2}
\end{prop} 

\noindent
{\em Proof.} Since  $\pi_i(v_j) = \pi_i(v_k)$, we have 
$v_j - v_k = av_i$ for some  $a \in \Z$. Obviously,  
$a \neq 0$. Assume that $a > 0$. Then the point $v_j= v_k + av_i$ belongs 
to the relative interior of the cone generated by $v_i$ and $v_k$. 
Therefore, by \ref{def.fano}(iii,iv), the set ${\cal P} := \{v_i, v_k\}$ must 
be a primitive collection. The corresponding primitive relation 
can not be  $v_i + v_k =0$, because otherwise $\pi_i(v_k)$ would be $0$. 
By \ref{coeff-prim}, the only possibility for the primitive relation 
$R({\cal P})$ is $v_i + v_k = v_l$. This implies that $\pi_i(v_l)=
\pi_i(v_k) = v$. By \ref{def.fano}(iv), $\pi_i^{-1}(v)$ can't contain more
than $2$ vertices of $P$. Hence, $v_l = v_j$ and $a =1$. So we come 
to the situation (i). Analogously, if $a<0$, then one obtains the second 
case (ii). 
\hfill $\Box$ 

\begin{dfn} 
{\rm Let $v \in P_i= \pi_i(P)$ be a nonzero double point ($0 \neq 
v = \pi_i(v_j) = \pi_i(v_k))$ . If 
$v_i + v_k = v_j$ is a primitive relation, then we call the vertex 
$v_j \in P$ a {\bf $\pi_i$-link of $v_k$}. 
If $v_i + v_j = v_k$ is a primitive relation, then we call the vertex 
$v_k \in P$ a {\bf $\pi_i$-link of $v_j$}.}
\end{dfn}

\begin{prop}
Let  $v_k$ be a vertex in ${\cal V}(P) \setminus \{v_i, v_j \}$. Then 
$\pi_i(v_k)$ is a double point of $P_i$ if and only if exactly one 
of two sets $\{ v_i, v_k\}$ and $\{v_j,v_k\}$ is a primitive 
collection.
\label{cent-sym2}  
\end{prop} 

\noindent
{\em Proof.} Assume that $\pi_i(v_k)$ is a double point of $P_i$, i.e., 
there exists a vertex $v_l \in  {\cal V}(P) \setminus \{v_i, v_j, v_k \}$
such that $\pi_i(v_l) = \pi_i(v_k)$. By \ref{double-2}, exactly 
one of the two subsets $\{v_i, v_k\}$ and $\{v_i, v_l\}$ is a primitive 
collection. Analogously, by \ref{double-2}, exactly 
one of the two subsets $\{v_j, v_k\}$ and $\{v_j, v_l\}$ is a primitive 
collection. If $\{v_j, v_l\}$ and  $\{v_i, v_l\}$ were primitive collections, 
then, by  \ref{double-2}, we would have $v_j + v_l = v_k$ and 
$v_i + v_l = v_k$. This implies a  contradiction to $v_i \neq v_j$. 
 If both $\{v_j, v_k\}$ and  $\{v_i, v_k\}$ were primitive collections, 
then, by  \ref{double-2}, we would have $v_j + v_k = v_l$ and 
$v_i + v_k = v_l$. This contradicts $v_i \neq v_j$. 

Now assume that  $\pi_i(v_k)$ is not a double point of $P_i$.
If e.g.  
$\{ v_k, v_i \}$ is a primitive collection, then the corresponding 
primitive relation can't be of the form $v_k + v_i =0$ (otherwise we would 
have $v_i = v_j$). Hence, the only possiblity for the primitive collection 
is $v_k + v_i = v_l$ where  $v_l \in {\cal V}(P) \setminus \{v_i, v_j, v_k 
\}$.  This implies that $\pi_i(v_k) = \pi_i(v_l)$, i.e., 
$\pi_i(v_k)$ is a double point of $P_i$. Contradiction. 

\hfill $\Box$

\subsection{Toric Fano $3$-folds}

The purpose of this section is to illustrate for  toric Fano $3$-folds 
the method which will be the main tool  for our  
classification of $4$-dimensional Fano manifolds. Our purpose 
is the following theorem which was proved in \cite{bat1} and \cite{wat}: 

\begin{theo}  
There exist exactly $18$ different types of   
$3$-dimensional smooth toric Fano varieties. The maximum of the Picard 
number of such varieties is equal to $5$. 
\end{theo} 

This statement can be equivalently reformulated as follows: 

\begin{theo}  
There exist exactly 
$18$ different types of $3$-dimensional Fano polyhedra $P$ up to isomorphism. 
The maximal number of vertices of $P$ is $8$:

\medskip

\begin{center}
\begin{tabular}{|c|c|c|c|c|c|c|c|c|} \hline
\mbox{\em the number of vertices \/} & 4 & 5 & 6 & 
7 & 8  \\ \hline 
\mbox{\em the number of polyhedra\/} & 1 & 4 & 
\makebox{7} & \makebox{4} & \makebox{2}  \\ \hline 
\end{tabular}
\end{center}
\end{theo}
\medskip

\begin{prop} 
Assume that a $3$-dimensional Fano polyhedron satisfies the property:
for every primitive collection $\{v_i, v_j \} \subset {\cal V}(P)$, the 
corresponding primitive relation is $v_i + v_j = 0$. Then the number $n = 
n(P)$ of vertices of $P$ is not greater than $6$. 
\label{sym.3}
\end{prop}  

\noindent 
{\em Proof .} Let $k$ be the number of primitive collections 
in ${\cal V}(P)$ consisting of $2$ vertices. Since $f_1(P) = 3n -6$,  
we have 
\[ k=  { n \choose 2 } - (3n - 6) = \frac{(n-3)(n-4)}{2}. \]
On the other hand, our assumption on  $P$ implies that any two different 
primitive collections consisting of $2$ vertices can not have common 
elements, i.e.,  $2k \leq n $. 
Therefore, $n^2 - 8n + 12 \leq 0$. This implies $n \leq 6$. 
\hfill $\Box$

\begin{prop} Any $3$-dimensional Fano polyhedron has 
at most $8$ vertices. 
\end{prop}

\noindent 
{\em Proof.}  Using  \ref{sym.3}, we can assume that there exists  
a primitive relation of the type 
\[ v_1 + v_2 - v_3 = 0. \]
By \ref{except}, $D_3$ is a ruled toric surface, i.e., $[v_3, v_1]$ and 
$[v_3,v_2]$ are faces of $P$ and there are at most  $2$ more vertices of 
$P$, e.g.,  $v_4$ and $v_5$, which are connected by edges $[v_3,v_4]$ 
and $[v_3,v_5]$ with $v_3$. 
By \ref{commun}, among $\{ \pi_3(v_1), 
\pi_3(v_2), \pi_3(v_4), \pi_3(v_5) \} \subset P_3$ there exist 
at most two double points. Thus, there might be 
at most three vertices of $P$ which are not joined with $v_3$ by an 
edge: one 
centrally symmetric vertex $v_6 = -v_3$, and  two vertices $v_7$ and 
$v_8$ such that $\pi_3(v_7)$ and $\pi_3(v_8)$ are double 
vertices of $P_3$ (the latter implies that 
$\{ v_3, v_7 \}$ and $\{ v_3, v_8 \}$  are 
primitive collections of degree  $1$). 
\hfill $\Box$

\begin{prop}
There exist exactly $4$ different $3$-dimensional Fano polyhedra 
with $5$ vertices.
\end{prop}

\noindent 
{\em Proof.}
It is known that the combinatorial type  of $3$-dimensional 
simplicial polyhedron  having $5$ vertices is unique: it is 
defined by two primitive collections $\{v_1,v_2\}$ and 
$\{v_3, v_4, v_5\}$.  Let  ${\cal B}$ 
denotes this combinatorial type. Then, up to change of indices 
in the numeration of vertices,  the corresponding 
primitive relations are      
\medskip

\begin{center}
\begin{tabular}{|c|c|c|c|c|} \hline
\FP & ${\cal B}_1$ & ${\cal B}_2$ & ${\cal B}_3$ & ${\cal B}_4$  \\ \hline
$v_1 + v_2  =$ & $$0$$ & $$0$$ & $v_3$ & $0$  \\ \hline 
$v_3 + v_4 + v_5 = $ & $2v_1$ & $v_1$ & $0$ & $0$ \\ \hline
\end{tabular}
\end{center}
\hfill $\Box$

\begin{prop} 
There exist exactly $7$ different $3$-dimensional Fano polyhedra 
with 6 vertices. 
\end{prop}

\noindent 
{\em Proof.}
It is easy to show that 
there exist exactly two different combinatorial types ${\cal C}$ and 
${\cal D}$ of $3$-dimensional Fano polyhedra   having $6$ vertices. 

Using explicit decription of arbitrary 
smooth projective toric varieties with the Picard number $3$ \cite{bat5} 
together with \ref{crit}, 
we obtain that, up to renumeration of vertices of $P$, the primitive 
relations fare 

1) the type ${\cal C}$:  $v_1 + v_2 = 0$ and 
\medskip

\begin{center}
\begin{tabular}{|c|c|c|c|c|c|} \hline
\FP & ${\cal C}_1$ & ${\cal C}_2$ & ${\cal C}_3$ & ${\cal C}_4$ & 
${\cal C}_5$  \\ \hline
$v_3 + v_5  =$ & $v_1$ & $v_1$  & $0$ & $0$ & $v_1$ \\ \hline
$v_4 + v_6 =$ & $v_1$ & $v_3$ & $0$ & $v_3$ & $v_2$ \\ \hline
\end{tabular}
\end{center}
\medskip

2) for the type ${\cal D}$:  
$v_3 + v_6 = 0$, $v_4 + v_6 = v_5$, $v_3 + v_5 = v_4$, and 
\medskip
\begin{center}
\begin{tabular}{|c|c|c|} \hline
\FP & ${\cal D}_1$ & ${\cal D}_2$  \\ \hline
$v_1 + v_2 + v_4=$ & $2v_3$ & $v_3$  \\ \hline
$v_1 + v_2 + v_5 =$ & $v_3$ & $0$  \\ \hline
\end{tabular}
\end{center}
\hfill $\Box$
\medskip

\begin{prop} 
Assume that a $3$-dimensional Fano polyhedron $P$ contains two 
centrally symmetric vertices $v_i, v_j \in {\cal V}(P)$: $v_i + v_j =0$. 
Then $P_i = \pi_i(P)$ is a $2$-dimensional Fano polyhedron. 
Moreoever, if $v + v' = v''$ holds for some vertices $v, v', v'' \in 
{\cal V}(P_i)$, then no one of two vertices $v$ and $v'$ is a double 
point of $P_i$.
\label{cent-3dim} 
\end{prop} 

\noindent
{\em Proof.} By \ref{cor-gamma}, if $P_i$ is not a Fano polyhedron, then 
there exists a face $\Gamma$ of $P_i$ containing exactly $3$ lattice points. 
The only possible linear relation (obtained form some primitive 
relation of degree $1$) satisfied by these $3$ lattice points 
could be 
$$\pi_i(v_{i_1}) + \pi_i(v_{i_2}) = 2\pi_i(v_{j_1}).$$ 
This relation contradicts \ref{cent-sym1}. 

Let $v, v', v''$ be vertices in   
${\cal V}(P_i)$ such that $v = \pi_i(v_k)$, $v'=  \pi_i(v_l)$,  
$v''=  \pi_i(v_m)$, and  $v + v' = v''$. Assume, for instance, that 
$v$ is a double point, i.e., there exists a vertex $v_s \in {\cal V}(P) 
\setminus \{v_k\}$ such that $\pi_i(v_s) =\pi_i(v_k) = v$. 
The linear relation  $v + v' = v''$ implies that at least one 
of the followg two equalities holds: 
\[ v_k + v_l = v_m + av_i, \;\; a \in \Z_{\geq 0}, \] 
\[  v_k + v_l = v_m + bv_j, \;\; b \in \Z_{\geq 0}. \]
Assume, for instance, that the first equality holds. We want to show 
that in this case 
$\{ v_k, v_l \} \subset {\cal V}(P)$ is a primitive collection. 
If it were  not the case, then we would have 
$a \neq 0$ and $[v_m, v_i]$ couldn't be  an edge
of $P$ (otherwise two cones $\sigma([v_k, v_l])$ and $\sigma([v_m,v_i])$ 
would have a common point in their relative interior). So, $\{ v_m,v_i \}$ 
is a primitive collection. On the other hand, the single  possibility 
for the corresponding primitive relation is $v_m + v_i = v_t$, i.e., 
$v''$ is another double point. This contradicts \ref{commun}. Hence, 
 $\{ v_k, v_l \}$ must be a primitive collection. 
The single possibility for corresponding primitive relation 
is  $v_k + v_l = v_m$. On the other hand, the same  arguments 
applied to the set of vertices $\{v_s, v_l, v_m\}$ instead of 
$\{v_k, v_l, v_m\}$ show that  $v_s + v_l = v_m$ must be 
a primitive relation too. Both above relations imply $v_k = v_s$. 

\hfill $\Box$

\begin{prop}  
There exist  exactly $2$  different $3$-dimensional Fano polyhedra 
with $8$ vertices. 
\end{prop}

\noindent 
{\em Proof.}
By \ref{sym.3}, there exists a primitive relation of the type $v_1 + v_3 = 
v_2$.  By \ref{ex-ray2}(i), $[v_1,v_2]$ and $[v_3,v_2]$ 
are edges of $P$. It follows from \ref{ex-ray2}(iii) 
that the divisor $D_2$ has 
Picard number $2$.  i.e., that $v_2$ is joined by edges of $P$ with exactly 
two more vertices, e.g.,  $v_7$ and $v_8$. By \ref{conv-h}, $P_2 = \pi_2(P)$ 
contains exactly $5$ lattice points. By \ref{commun}, at most $3$ lattice 
points in $P$ could be double points. Since $P$ contains $8$ points, $P_2$ 
must contain exactly $3$ double points. By \ref{commun}, these 
double points could be either $\{0, \pi_2(v_1), \pi_2(v_3) \}$, or 
 $\{0, \pi_2(v_7), \pi_2(v_8) \}$, i.e., 
there always exists a centrally symmetric to $v_2$ vertex $v_5 = - v_2$ and 
there always exist two vertices $v_4$, $v_6$ such that 
$\{ v_2, v_4 \}$ and $\{ v_2, v_6 \}$ are 
primitive collections. 
By \ref{cent-3dim}, $P_2$ must be a  Fano polygon and 
two double vertices of $P_2$ must be centrally symmetric.  
Without loss of generality we can assume that 
$v_1$ and $v_3$ are $\pi_2$-links of 
$v_4$ and $v_6$. Therefore the set   $\{ v_1, v_2, v_3, v_4, v_5, v_6 \}$ 
is  contained in a $2$-dimensional vector 
subspace. These  $6$ vertices of $P$ lying in this subspace form a 
$2$-dimensional Fano polyhedron with $6$ vertices and 
the combinatorial type of $P$ is determined uniquely (we denote it 
by ${\cal F}$). Now one sees that the primitive relations are 
$v_3 + v_6 = 0$, $v_1 + v_3 = v_2$, $v_2 + v_6 = v_1$, 
$v_2 + v_4 = v_3$, $v_2 + v_5 = 0$, $v_1 + v_4 = 0$, 
$v_4 + v_6 = v_5$, $v_1 + v_5 = v_6$, $v_3 + v_5 = v_4$ and  

\begin{center}
\begin{tabular}{|c|c|c|} \hline
\FP & ${\cal F}_1$ & ${\cal F}_2$   \\ \hline
$v_7 + v_8 =$ & $0$ & $v_1$  \\ \hline
\end{tabular}
\end{center}

\hfill $\Box$
\medskip

\begin{prop} 
There exist exactly $4$  different $3$-dimensional Fano polyhedra 
with $7$ vertices. 
\end{prop}

\noindent 
{\em Proof.}
By \ref{sym.3}, there exists a primitive relation of the type $v_1 + v_3 = 
v_2$. By \ref{ex-ray2}(i), $[v_1,v_2]$ and $[v_3,v_2]$ 
are edges of $P$. It follows from \ref{ex-ray2}(iii) 
that the divisor $D_2$ has 
Picard number $2$.  i.e., that $v_2$ is joined by edges of $P$ with exactly 
two more vertices, e.g.,  $v_6$ and $v_7$. By \ref{commun}, we have 
the following $3$ cases: 

{\sc Case 1.} One of two vertices $v_4$ or $v_5$ is centrally symmetric 
to $v_2$, e.g.,  $v_2 + v_4 = 0$. Then we use \ref{cent-3dim} to  
obtain that $P_2$ is a Fano polygon whose single double vertex 
$\pi_2(v_5)$ belongs to  
a pair of  centrally symmetric vertices of $P_2$.  
So,  one of two 
 sets $\{ v_1, v_2, v_3, v_4, v_5 \}$ or  
$\{  v_2, v_4, v_5, v_6, v_7 \}$ must be contained in a $2$-dimensional vector 
subspace. Therefore, vertices of $P$ lying in this subspace form a 
Fano polygon $P_2$  with $5$ vertices. We assume that 
vertices $\{ v_1, v_2, v_3, v_4, v_5 \}$ form such a polygon. Then 
the combinatorial type of $P$ is determined uniquely (we denote it 
by ${\cal E}$).  This shows that the 
 primitive relations up to a numeration of vertices 
are 
 $v_2 + v_4 = 0$, $v_3 + v_5 = 0$, $v_1  + v_3 = v_2$, 
$v_2 + v_5 = v_1$, $v_1 + v_4 = v_5$  and 

\begin{center}
\begin{tabular}{|c|c|c|c|c|} \hline
\FP & ${\cal E}_1$ & ${\cal E}_2$ & ${\cal E}_3$ & ${\cal E}_4$  \\ \hline
$v_6 + v_7 =$ & $v_1$ & $v_2$ &  $0$ & $v_3$ \\ \hline 
\end{tabular}
\end{center}

{\sc Case 2.} $\pi_2(v_1)$ and  $\pi_2(v_3)$ are double vertices 
of $P_2$, i.e., we have two primitive relations 
$\{ v_1, v_2, v_3, v_4, v_5 \}$
$v_2 + v_4 = v_1$ and $v_2 + v_5 = v_2$.  Since   
$\pi_2(v_1)$ and  $\pi_2(v_3)$
are centrally symmetric, the   set $\{ v_1, v_2, v_3, v_4, v_5 \}$
is contained in a $2$-dimensional vector 
subspace and 
at least one of two 
vertices $v_1$ and  $v_2$ must possess a centrally symmetric one, we come  
to the situation considered in Case 1. 

{\sc Case 3.} $\pi_2(v_6)$ and  $\pi_2(v_7)$ are double vertices 
of $P_2$, i.e., we have two primitive relations
$v_2 + v_4 = v_6$ and $v_2 + v_5 = v_7$. If  $\pi_2(v_6)$ and  $\pi_2(v_7)$
are centrally symmetric, then the set 
$\{  v_2, v_4, v_5, v_6, v_7 \}$ is contained in a $2$-dimensional vector 
subspace and at least one of two 
vertices $v_6$ and  $v_7$ must possess a centrally symmetric one, i.e., 
we come  to the situation considered in Case 1. It remains to exclude
the situation, when  $\pi_2(v_6)$ and  $\pi_2(v_7)$ are not centrally 
symmetric vertices of $P_2$. Indeed, by \ref{ex-ray2}(i), both 
vertices $v_6$ and $v_7$ have valence $4$. Moreover, 
$\{v_6, v_5\}$, $\{v_6,v_7\}$ 
  and $\{v_4,v_7\}$ are primitive collections. Since
$$\pi_2(v_6) + \pi_2(v_5) = \pi_2(v_6) + \pi_2(v_7) = 
\pi_2(v_4) + \pi_2(v_7) = v \neq 0$$
where $v$ is not a double vertex of $P_2$, we obtain the single
possibility for the corresponding primitive relations: 
\[ v_6 + v_5 = v_6 + v_7 = v_4 + v_7 = v_i, \]
where $v_i$  is a vertex of $P$. The last equalities are 
impossible, since they imply $v_5=v_7$ and $v_6 = v_4$. 

\hfill $\Box$
\medskip

\begin{rem}
{\rm In the table below we give the  list of 
3-dimensional toric Fano varieties $V=V(P)$. We denote by $S_i$ the 
Del Pezzo surface obtained by the blow up of $i$ points on ${\P}^2$. 
The number  $a(V)$ denotes the dimension of the group 
${\rm Aut}\; (V)$ of biregular authomorphisms of $V$. }
\end{rem}
\medskip

\begin{center}
\begin{tabular}{|c|c|c|c|c|l|} \hline 
    &      &     &     &     &      \\
 $n^0$ &  $c_1^3$  & $b_2$ & $h^0$  & $a(V)$ &
the type of Fano polytope $P$
\\    &      &     &     &     &  and the geometry of $V$    
\\ \hline 
  1 &  64 & 1 &  35   & 15 &  ${\P}^3$
\\    &      &     &     &     &      
\\ \hline
  2 &  62 & 2 &  34  &  15  &  ${\cal B}_1, \; {\P}_{{\P}^2} 
({\cal O} \oplus {\cal O}(2))$
\\    &      &     &     &     &      
\\ \hline
  3 &  56 & 2 &  31  &  12  &  ${\cal B}_2, \; {\P}_{{\P}^2} 
({\cal O} \oplus {\cal O}(1))$
\\    &      &     &     &     &      
\\ \hline
  4 &  54 & 2 &  30  &  11  &  ${\cal B}_3, \; {\P}_{{\P}^1}
({\cal O} \oplus {\cal O} \oplus {\cal O}(1))$
\\    &      &     &     &     &      
\\ \hline
  5 &  54 & 2 &  30  &  11  &  ${\cal B}_4 , \; {\P}^2 \times {\P}^1$
\\    &      &     &     &     &      
\\ \hline
  6 &  52 & 3 &  29  &  11 & ${\cal C}_1$, 
${\P}_{{\P}^1 \times {\P}^1} ({\cal O} \oplus 
{\cal O}(1,1))$
\\    &      &     &     &     &      
$\;\;\;\;$  
\\ \hline
  7 &  50 & 3 &  28  & 10  & ${\cal C}_2 , \; {\P}_{S_1} 
  ({\cal O} \oplus {\cal O}(l) $, 
 where $l^2 = 1$ on $S_1$    
\\    &      &     &     &     & 
\\ \hline
  8 &  48 & 3 &  27  &  9  & ${\cal C}_3 , \;\npo \times \npo \times \npo$
\\    &      &     &     &     &      
\\ \hline
 9 &  48 & 3 &  27  &  9   & ${\cal C}_4 , \;S_1  \times {\P}^1$
\\    &      &     &     &     &     
\\ \hline
 10 &  44 & 3 &  25  &  7  &  ${\cal C}_5$, 
${\P}_{\npo \times \npo} ({\cal O} \oplus 
{\cal O}(1,-1))$
\\    &      &     &     &     & 
$\;\;\;\;$   
\\ \hline
11 &  50 & 3 &  28  &  10  & ${\cal D}_1$, the blow up of $\npo$ on $n^03$
\\    &      &     &     &     &     
\\ \hline

12 &  46 & 3 &  26  &  8   & ${\cal D}_2$,  the blow up of $\npo$  on $n^05$
\\    &      &     &     &     &     
\\ \hline
13 &  46 & 4 &  26  &  9  &  ${\cal E}_1$, $S_2$-bundle over $\npo$ 
\\    &      &     &     &     &    
\\ \hline
14 &  44 & 4 &  25  &  8  &  ${\cal E}_2$, $S_2$-bundle over $\npo$ 
\\    &      &     &     &     &     
\\ \hline
15 &  42 & 4 &  24  &  7  &  ${\cal E}_3$, $S_2 \times \npo$
\\    &      &     &     &     &    
\\ \hline
16 &  40 & 4 &  23  &  6  &  ${\cal E}_4$, $S_2$-bundle over $\npo$
\\    &      &     &     &     &      
\\ \hline
17 &  36 & 5 &  21  &  5  &  ${\cal F}_1$, $S_3 \times \npo$
\\    &      &     &     &     &     
\\ \hline
18 &  36 & 5 &  21 &  5  &  ${\cal F}_2$, $S_3$-bundle over $\npo$ 
\\    &      &     &     &     &     
\\ \hline  
\end{tabular}
\end{center}
\bigskip

\section{Classification of $4$-dimensional Fano polyhedra}

\subsection{Fano $4$-polyhedra with $\leq 7$ vertices}

There there exists a unique 
toric $4$-fold, $4$-dimensional projective space ${\P}^4$, having the 
Picard number $1$. The corresponding Fano $4$-polyhedron is a simplex.

\begin{prop}
There exist exactly $9$ different $4$-dimensional Fano polyhedra with 
$6$ vertices having two possible combinatorial types $B$ and $C$. These 
$4$-polyhedra are defined by the following primitive relations

\begin{center}
\begin{tabular}{|c|c|c|c|c|c|} \hline
\FP & $B_1$ & $B_2$ & $B_3$ & $B_4$ & $B_5$ \\ \hline
$v_1 + v_2 + v_3 + v_4 =$ & $3v_5$ & $2v_5$ & $v_5$ & $0$ & $0$ \\ \hline 
$v_5 + v_6 = $ & $0$ & $0$ & $0$ & $0$ & $v_1$ \\ \hline
\end{tabular}
\end{center}

\medskip
\begin{center}
\begin{tabular}{|c|c|c|c|c|} \hline
\FP & $C_1$ & $C_2$ & $C_3$ & $C_4$ \\ \hline
$v_1 + v_2 + v_3 =$ & $0$ & $0$ & $0$ & $0$ \\ \hline
$v_4 + v_5 + v_6 =$ & $2v_1$ & $v_1$ & $v_1 + v_2$ & $0$ \\ \hline
\end{tabular}
\end{center}
\end{prop}

\noindent 
{\em Proof.} The statement immediately follows from the general 
description of toric Fano manifolds with the Picard number $2$ (see  
\cite{kl}) and from \ref{crit}. 
\hfill $\Box$

\begin{prop}
\label{poly7}
There exist exactly $28$ different $4$-dimensional Fano polyhedra with 
$7$ vertices having three  possible combinatorial types $D$, $E$, and $G$. 
These $4$-polyhedra are defined by the following primitive relations

{\rm (i)} Type $D:$

\begin{center}
\begin{tabular}{|c|c|c|c|c|c|c|c|} \hline
\FP & $D_1$ & $D_2$ & $D_3$ & $D_4$ & $D_5$ & $D_6$ & $D_7$ \\ \hline
$v_1 + v_2 + v_3=$ & $2v_6$ & $2v_4$ & $v_4 + v_6$ & $2v_6$ & $2v_4$ & $v_6$ 
& $0$ \\ \hline
$v_4 + v_5 =$ & $v_6$ & $v_6$ & $v_6$ & $v_1$ & $0$ & $v_6$ & $v_1$ \\ \hline
$v_6 + v_7 =$ & $0$ & $0$ & $0$ & $0$ & $0$ & $0$ & $v_1$ \\ \hline
\end{tabular}
\end{center}
\medskip

\begin{center}
\begin{tabular}{|c|c|c|c|c|c|c|} \hline
\FP & $D_8$ & $D_9$ & $D_{10}$ & $D_{11}$ & $D_{12}$ & $D_{13}$  \\ \hline
$v_1 + v_2 + v_3=$ & $v_4$ & $v_4 + v_6$ & $v_6$ & $0$ & $v_4$ & $0$  
\\ \hline
$v_4 + v_5 =$ & $v_6$ & $0$ & $v_1$ & $v_1$ & $0$ & $0$  \\ \hline
$v_6 + v_7 =$ & $0$ & $0$ & $0$ & $v_4$ & $0$ & $0$  \\ \hline
\end{tabular}
\end{center}
\medskip

\begin{center}
\begin{tabular}{|c|c|c|c|c|c|c|} \hline
\FP & $D_{14}$ & $D_{15}$ & $D_{16}$ & $D_{17}$ & $D_{18}$ & $D_{19}$  
\\ \hline
$v_1 + v_2 + v_3=$ & $0$ & $0$ & $v_4 + v_7$ & $0$ & $2v_7$ & $v_7$  
\\ \hline
$v_4 + v_5 =$ & $0$ & $0$ & $v_6$ & $v_1$ & $v_6$ & $v_6$  \\ \hline
$v_6 + v_7 =$ & $v_1$ & $v_4$ & $0$ & $v_2$ & $0$ & $0$  \\ \hline
\end{tabular}
\end{center}

{\rm (ii)} Type  $E:$  
$v_1 + v_7 = 0$, $v_1 + v_2 = v_6$, $v_6  + v_7 = v_2$  and
 
\begin{center}
\begin{tabular}{|c|c|c|c|} \hline
\FP & $E_1$ & $E_2$ & $E_3$  \\ \hline
$v_2 + v_3 + v_4 + v_5 =$ & $2v_1$ & $v_1$ &  $0$  \\ \hline 
$v_3 + v_4 + v_5 + v_6 =$ & $3v_1$ & $2v_1$ & $v_1$  \\ \hline
\end{tabular}
\end{center}
\medskip

{\rm (iii)} Type $G:$ 

\begin{center}
\begin{tabular}{|c|c|c|c|c|c|c|} \hline
\FP & $G_1$ & $G_2$ & $G_3$ & $G_4$ & $G_5$ & $G_6$  \\ \hline
$v_1 + v_7 =$ & $0$ & $v_4$ & $0$ & $v_4$ & $v_4$ & $v_4$  \\ \hline
$v_2 + v_3 + v_4 = $ & $v_1$ & $v_7$ & $0$  & $v_7$ & $v_7$ & $v_7$ 
 \\ \hline
$v_4 + v_5 + v_6 = $ & $2v_1$ & $2v_1$ & $v_1$ & $v_1 + v_2$ & $0$ & $v_1$ 
\\ \hline
$v_5 + v_6 + v_7=$ & $v_2 + v_3$ & $v_1$ & $v_2 + v_3$ & $v_2 $ &  
 $v_2 + v_3$ & $0$ \\ \hline
$v_1 + v_2 + v_3=$ & $v_5 + v_6$ & $0$ & $v_5 + v_6$ & $0$ & $0$ & $0$ 
\\ \hline
\end{tabular}
\end{center}
\end{prop}

\noindent 
{\em Proof.}  There are $4$ different posssible combinatorial types of 
simplicial fans with $7$ vertices defining smooth projective toric $4$-folds. 
These combinatorial types are described by the following Gale diagrams:

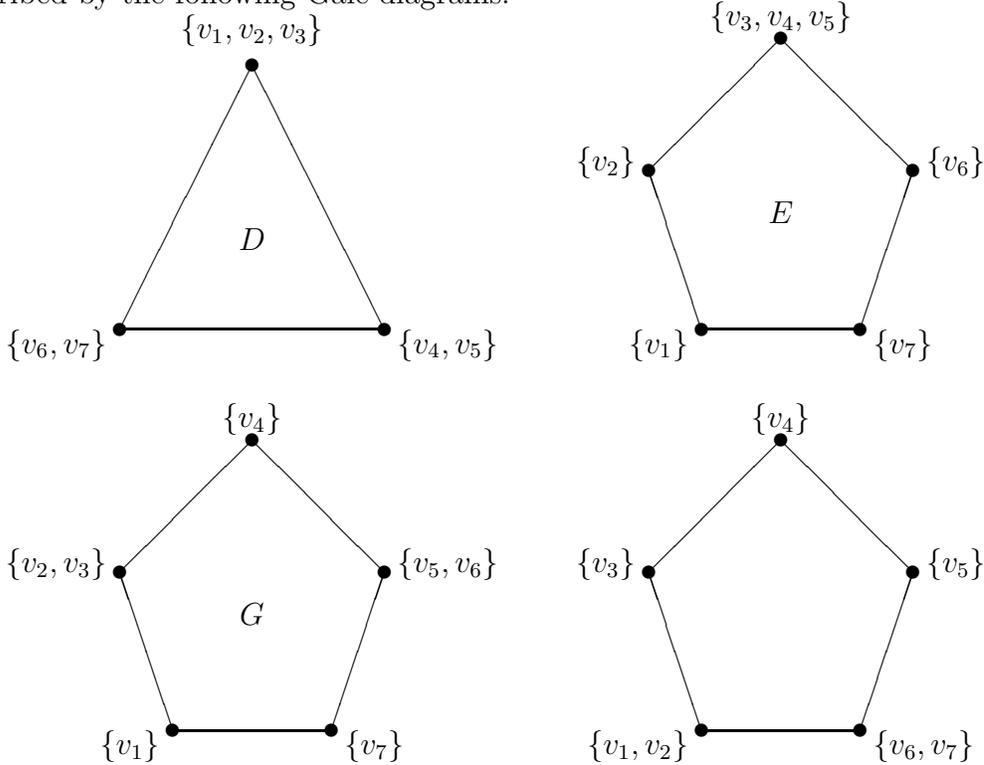
\begin{picture}(500,140)

\put(50,20){\circle*{5}} 
\put(150,20){\circle*{5}}
\put(100,120){\circle*{5}}

\put(100,130){\makebox(0,0)[b]{$\{v_1,v_2,v_3\}$}}
\put(45,20){\makebox(0,0)[tr]{$\{v_6,v_7\}$}}
\put(155,20){\makebox(0,0)[tl]{$\{v_4,v_5\}$}}
\put(50,20){\line(1,0){100}}
\put(50,20){\line(1,2){50}}
\put(150,20){\line(-1,2){50}}
\put(100,50){\makebox(0,0)[b]{$D$}}

\put(270,20){\circle*{5}} 
\put(330,20){\circle*{5}}
\put(250,80){\circle*{5}}
\put(350,80){\circle*{5}} 
\put(300,130){\circle*{5}}

\put(270,20){\line(1,0){60}}
\put(270,20){\line(-1,3){20}}
\put(330,20){\line(1,3){20}}
\put(300,130){\line(-1,-1){50}}
\put(300,130){\line(1,-1){50}}

\put(300,135){\makebox(0,0)[b]{$\{v_3,v_4,v_5\}$}}
\put(265,20){\makebox(0,0)[tr]{$\{v_1\}$}}
\put(335,20){\makebox(0,0)[tl]{$\{v_7\}$}}
\put(245,80){\makebox(0,0)[br]{$\{v_2\}$}}
\put(355,80){\makebox(0,0)[bl]{$\{v_6\}$}}

\put(300,60){\makebox(0,0)[b]{$E$}}

\end{picture}

\bigskip

\begin{picture}(500,140)
\put(70,20){\circle*{5}} 
\put(130,20){\circle*{5}}
\put(50,80){\circle*{5}}
\put(150,80){\circle*{5}} 
\put(100,130){\circle*{5}}

\put(70,20){\line(1,0){60}}
\put(70,20){\line(-1,3){20}}
\put(130,20){\line(1,3){20}}
\put(100,130){\line(-1,-1){50}}
\put(100,130){\line(1,-1){50}}

\put(100,135){\makebox(0,0)[b]{$\{v_4\}$}}
\put(65,20){\makebox(0,0)[tr]{$\{v_1\}$}}
\put(135,20){\makebox(0,0)[tl]{$\{v_7\}$}}
\put(45,80){\makebox(0,0)[br]{$\{v_2,v_3\}$}}
\put(155,80){\makebox(0,0)[bl]{$\{v_5,v_6\}$}}

\put(100,60){\makebox(0,0)[b]{$G$}}

\put(270,20){\circle*{5}} 
\put(330,20){\circle*{5}}
\put(250,80){\circle*{5}}
\put(350,80){\circle*{5}} 
\put(300,130){\circle*{5}}

\put(270,20){\line(1,0){60}}
\put(270,20){\line(-1,3){20}}
\put(330,20){\line(1,3){20}}
\put(300,130){\line(-1,-1){50}}
\put(300,130){\line(1,-1){50}}

\put(300,135){\makebox(0,0)[b]{$\{v_4\}$}}
\put(265,20){\makebox(0,0)[tr]{$\{v_1,v_2\}$}}
\put(335,20){\makebox(0,0)[tl]{$\{v_6,v_7\}$}}
\put(245,80){\makebox(0,0)[br]{$\{v_3\}$}}
\put(355,80){\makebox(0,0)[bl]{$\{v_5\}$}}
\end{picture}

Using the general description 
of toric Fano manifolds with the Picard number $3$ \cite{bat5} and 
from \ref{crit}, it is easy to show that
only first $3$ combinatorial types admit a realization by a Fano 
$4$-polyhedron.  
\hfill $\Box$

\subsection{Fano $4$-polyhedra having a vertex of valence $5$}

\begin{prop} Let $P$ be a Fano $4$-polyhedron, $v_0$ an 
a vertex of $P$, $P_0 = \pi_0(P)$ the projection of the polyhedron 
$P$ in ${\R}^3 / {\R}\langle v_0 \rangle$. 
Then there exist  only $4$ possibilities 
for $2$-dimensional faces of the polyhedron   $P_0$: 

\bigskip
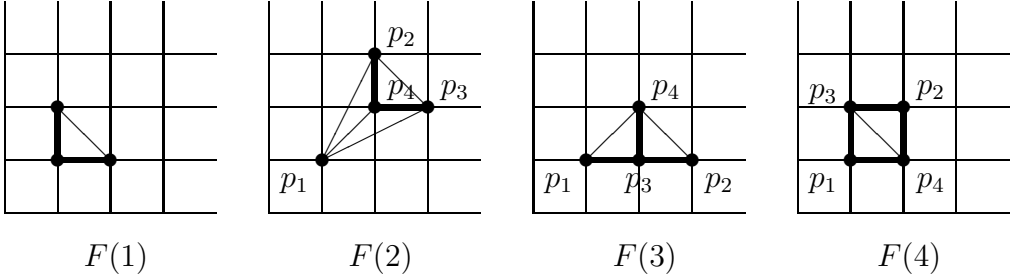
\begin{picture}(500,120)
\multiput(0,20)(20,0){4}%
{\line(0,1){80}}
\multiput(0,20)(0,20){4}%
{\line(1,0){80}}

\multiput(100,20)(20,0){4}%
{\line(0,1){80}}
\multiput(100,20)(0,20){4}%
{\line(1,0){80}}

\multiput(200,20)(20,0){4}%
{\line(0,1){80}}
\multiput(200,20)(0,20){4}%
{\line(1,0){80}}

\multiput(300,20)(20,0){4}%
{\line(0,1){80}}
\multiput(300,20)(0,20){4}%
{\line(1,0){80}}

\linethickness{0.7mm}
\put(20,40){\line(1,0){20}}
\put(20,60){\line(1,-1){20}}
\put(20,40){\line(0,1){20}}

\put(120,40){\line(2,1){40}}
\put(120,40){\line(1,2){20}}
\put(140,80){\line(1,-1){20}}
\put(140,60){\line(1,0){20}}
\put(140,60){\line(0,1){20}}
\put(140,60){\line(-1,-1){20}}

\put(220,40){\line(1,0){40}}
\put(220,40){\line(1,1){20}}
\put(260,40){\line(-1,1){20}}
\put(240,60){\line(0,-1){20}}

\put(320,40){\line(1,0){20}}
\put(320,40){\line(0,1){20}}
\put(340,60){\line(-1,0){20}}
\put(340,60){\line(0,-1){20}}
\put(320,60){\line(1,-1){20}}

\put(20,40){\circle*{5}}
\put(20,60){\circle*{5}}
\put(40,40){\circle*{5}}

\put(120,40){\circle*{5}}
\put(140,60){\circle*{5}}
\put(140,80){\circle*{5}}
\put(160,60){\circle*{5}}
\put(115,35){\makebox(0,0)[tr]{$p_1$}}
\put(145,65){\makebox(0,0)[bl]{$p_4$}}  
\put(145,85){\makebox(0,0)[bl]{$p_2$}}  
\put(165,65){\makebox(0,0)[bl]{$p_3$}}  

\put(220,40){\circle*{5}}
\put(240,40){\circle*{5}}
\put(260,40){\circle*{5}}
\put(240,60){\circle*{5}}

\put(215,35){\makebox(0,0)[tr]{$p_1$}}
\put(245,65){\makebox(0,0)[bl]{$p_4$}}  
\put(240,35){\makebox(0,0)[t]{$p_3$}}  
\put(265,35){\makebox(0,0)[tl]{$p_2$}}

\put(320,40){\circle*{5}}
\put(340,40){\circle*{5}}
\put(320,60){\circle*{5}}
\put(340,60){\circle*{5}}

\put(315,35){\makebox(0,0)[tr]{$p_1$}}
\put(345,65){\makebox(0,0)[bl]{$p_2$}}  
\put(345,35){\makebox(0,0)[tl]{$p_4$}}  
\put(315,65){\makebox(0,0)[rp]{$p_3$}}

\put(30,0){{$F(1)$}} 
\put(130,0){{$F(2)$}} 
\put(230,0){{$F(3)$}} 
\put(330,0){{$F(4)$}} 

\end{picture}

Moreover, if $P$ contains a centrally symmetric vertex $v_i = - v_0$, then 
only $F(1)$ and $F(4)$ can be $2$-dimensional faces of $P_0$. 
\label{3-dim.face}
\end{prop} 

\noindent 
{\em Proof. } We apply \ref{cor-gamma}(ii)  to the $4$-dimensional 
case. It is sufficient to find all possible primitive relations $R({\cal P})$ 
with $\Delta({\cal P}) =1$ giving  rise to faces $\Gamma \subset P_i$ 
which are diffrent from the standard triangle $F(1)$. 
Using \ref{ex-ray2}(i), we find the 
following  types of primitive relations:  
\[ v_0 + v_1 + v_2 + v_3 = 3 v_4\;\; \; \mbox{\rm  $\to F(2)$}, \;\;\;\;\;
v_0 + v_1 + v_2  = 2 v_3 \;\; \; \mbox{\rm  $\to F(3)$}, \] 
and 
\[ v_0 + v_1 + v_2  = v_3 +  v_4 \;\; \; 
\mbox{\rm  $\to F(4)$}. \]
By \ref{cent-sym1}, if there exists a vertex $v_i = -v_0$, then 
only the primitive relation 
$v_0 + v_1 + v_2  = v_3 +  v_4$ is possible. 
\hfill $\Box$

It turns out that the existence of two centrally symmetric vertices 
$v_i$ and $v_j$ 
of a Fano $4$-polyhedron $P$ allows completely 
describe the combinatorial structure 
of $P$ via its $\pi_i$-projection $P_i$: 

\begin{theo} 
Assume that $P$ contains two centrally symmetric vertices $v_i$ and $v_j$, 
i.e., $v_i + v_j = 0$ is a primitive relation. Then a set 
$\{v_{1}, v_{2},v_{3}, v_{4} \} \subset {\cal V}(P)$ generate a 
$4$-dimensional cone $\sigma(F)$ 
in $\Sigma(P)$ if and only if we have  one of the 
following situations: 

{\rm (i)} the set $\{v_{1},  v_{2},v_{3}, v_{4} \}$ contains $v_i$ 
$($or resp. $v_j$$)$, the other $3$ vertices are connected by eadges with 
$v_i$ $($or resp. with $v_j$),  and the convex hull of their  
$\pi_i$-projections is a $2$-dimensional simplicial face of $P_i$ of 
the type $F(1)$; 

{\rm (ii)} $\{v_{1},  v_{2},v_{3}, v_{4} \} 
\cap \{v_i, v_j \} = \emptyset$ and 
there exist two different vertices $v_k, v_l \in 
\{v_{1},  v_{2},v_{3}, v_{4} \}$ such that 
$\pi_i(v_k) = \pi_(v_l)$, and 
 the convex hull of  
$$\{\pi_i(v_{1}), \pi_i(v_{2}), \pi_i(v_{3}),  
\pi_i(v_{4}) \}$$ is 
a $2$-dimensional simplicial face $\Gamma$ of $P_i$ of the type $F(1)$; 

{\rm (iii)} the set $\{v_{1}, v_{2}, v_{3}, v_{4} \}$ contains $v_i$ 
$($or resp. $v_j$$)$,  the convex hull of the $\pi_i$-projections of  
the other $3$ vertices is a $2$-dimensional simplex $\Theta$ 
$($of the type $F(1))$  which is a half of a face $\Gamma \subset 
P_i$ of the type $F(4)$, and the cone over $\Theta$ belongs to 
$\Sigma_i(P)$ $($resp. to $\Sigma_j(P)$$)$;

{\rm (iv)}  $\{v_{1},  v_{2},v_{3}, v_{4} \} 
\cap \{v_i, v_j \} = \emptyset$, 
each vertex from  $\{v_{1},  v_{2}, v_{3}, v_{4} \}$ is
connected by an edge 
with both vertices $v_i, v_j$, and 
the convex hull of   
$$\{\pi_i(v_{1}),  \pi_i(v_{2}), \pi_i(v_{3}), \pi_i(v_{4})\}$$ 
is a $2$-dimensional face $\Gamma \subset P_i$ of the type $F(4)$, i.e., 
its  vertices satisfy a relation of the form 
\[ \pi_j(v_{k_1}) + \pi_j(v_{k_2}) = 
 \pi_j(v_{l_1}) +  \pi_j(v_{l_2}) \]
where $\{v_{k_1}, v_{k_2} \}$ and  $\{v_{l_1},  v_{l_2} \}$ 
are two disjoint subsets in  $\{v_{1},  v_{2},v_{3}, v_{4} \}$. 
\label{z-sym}
\end{theo}

\begin{coro} 
Let $P$ be a Fano $4$-polyhedron containing two centrally symmetric vertices 
$v_i, v_j  \in {\cal V}(P)$. Then the combinatorial type 
is uniquely determined
by the following data:

{\rm (i)} the combinatorial type of the 
$3$-dimensional fan $\Sigma_i(P)$ defining 
the toric divisor $D_i \subset V(P)$; 

{\rm (ii)} the convex polyhedron  $P_i = \pi_i(P)$; 

{\rm (iii)} the set $\delta(P_i) \subset {\cal V}(P_i)$ 
of all double vertices of $P_i$. 

In particular, the number $f_3(P)$ of $3$-dimensional faces 
of $P$ equals 
\[ 4 \alpha_i  + \beta_i -8  + \sum_{v \in \delta(P)} \gamma(v), \]
where $\alpha_i$ is the valence of $v_i$ in $P$, $\beta_i$ 
is the number of $2$-dimensional faces of $P_i$ of the type 
$F(4)$, and $\gamma(v)$ is the 
valence of a double vertex $v \in P_i$.    
\label{un-comb} 
\end{coro}

\begin{dfn} {\rm We use the notation $b, c_1, c_2$ and $d_2$ from 
\cite{grun} for the combinatorial types of some simplicial triangulations 
of the $2$-dimensional sphere $S^2$ with $5,6$ and $7$ vertices $n_i$.
There combinatorial types are described by the following  stereographic 
projections of $S^2$ from the vertex $n_1 = \infty$:

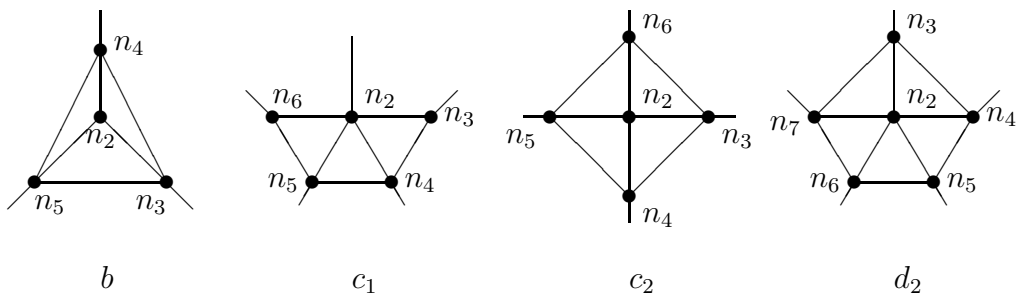
\begin{picture}(500,120)
\put(20,40){\circle*{5}} 
\put(70,40){\circle*{5}}
\put(45,65){\circle*{5}}
\put(45,90){\circle*{5}}
\put(20,35){\makebox(0,0)[tl]{$n_5$}} 
\put(70,35){\makebox(0,0)[tr]{$n_3$}} 
\put(45,60){\makebox(0,0)[t]{$n_2$}}
\put(50,90){\makebox(0,0)[bl]{$n_4$}}  

\put(20,40){\line(1,0){50}}
\put(20,40){\line(1,1){25}}
\put(20,40){\line(1,2){25}}
\put(20,40){\line(-1,-1){10}}
\put(70,40){\line(1,-1){10}}
\put(70,40){\line(-1,1){25}}
\put(70,40){\line(-1,2){25}}
\put(45,65){\line(0,1){40}}

\put(110,65){\circle*{5}}
\put(125,40){\circle*{5}}
\put(140,65){\circle*{5}}
\put(155,40){\circle*{5}}
\put(170,65){\circle*{5}}

\put(110,70){\makebox(0,0)[bl]{$n_6$}}
\put(120,40){\makebox(0,0)[r]{$n_5$}} 
\put(145,70){\makebox(0,0)[bl]{$n_2$}}
\put(160,40){\makebox(0,0)[l]{$n_4$}}  
\put(175,65){\makebox(0,0)[l]{$n_3$}}

\put(140,65){\line(-3,-5){20}}
\put(140,65){\line(3,-5){20}}
\put(140,65){\line(1,0){30}}
\put(140,65){\line(-1,0){30}}
\put(140,65){\line(0,1){30}}
\put(125,40){\line(1,0){30}}
\put(125,40){\line(-3,5){15}}
\put(155,40){\line(3,5){15}}
\put(110,65){\line(-1,1){10}}
\put(170,65){\line(1,1){10}}

\put(215,65){\circle*{5}}
\put(245,65){\circle*{5}}
\put(245,35){\circle*{5}}
\put(245,95){\circle*{5}}
\put(275,65){\circle*{5}}

\put(210,60){\makebox(0,0)[tr]{$n_5$}}
\put(250,70){\makebox(0,0)[bl]{$n_2$}}
\put(250,30){\makebox(0,0)[tl]{$n_4$}}
\put(280,60){\makebox(0,0)[tl]{$n_3$}}  
\put(250,100){\makebox(0,0)[l]{$n_6$}} 

\put(245,65){\line(0,-1){40}}
\put(245,65){\line(1,0){40}}
\put(245,65){\line(-1,0){40}}
\put(245,65){\line(0,1){40}}
\put(245,95){\line(1,-1){30}}
\put(245,95){\line(-1,-1){30}}

\put(215,65){\line(1,-1){30}}
\put(275,65){\line(-1,-1){30}}

\put(315,65){\circle*{5}}
\put(330,40){\circle*{5}}
\put(345,65){\circle*{5}}
\put(345,95){\circle*{5}}
\put(360,40){\circle*{5}}
\put(375,65){\circle*{5}}

\put(310,65){\makebox(0,0)[tr]{$n_7$}}
\put(325,40){\makebox(0,0)[r]{$n_6$}} 
\put(350,70){\makebox(0,0)[bl]{$n_2$}}
\put(365,40){\makebox(0,0)[l]{$n_5$}}  
\put(380,65){\makebox(0,0)[l]{$n_4$}}  
\put(350,100){\makebox(0,0)[l]{$n_3$}}

\put(345,65){\line(-3,-5){20}}
\put(345,65){\line(3,-5){20}}
\put(345,65){\line(1,0){30}}
\put(345,65){\line(-1,0){30}}
\put(345,65){\line(0,1){40}}
\put(345,95){\line(1,-1){30}}
\put(345,95){\line(-1,-1){30}}

\put(330,40){\line(1,0){30}}
\put(330,40){\line(-3,5){15}}
\put(360,40){\line(3,5){15}}
\put(315,65){\line(-1,1){10}}
\put(375,65){\line(1,1){10}}

\put(45,0){{$b$}} 
\put(140,0){{$c_1$}} 
\put(245,0){{$c_2$}} 
\put(345,0){{$d_2$}} 

\end{picture}

\label{stereo}
}
\end{dfn}

Using the restriction on the set of double points (see \ref{commun}), 
one obtains:  

\begin{theo} Let $P$ be a Fano $4$-polyhedron  with at least $8$ 
vertices and $v_0 \in P$ is  a vertex  of valence $5$. Assume that  
that vertices $n_j$ $(1 \leq j \leq 5)$ 
of the diagram \ref{stereo}$(b)$ 
correspond to vertices $v_j$ $(1 \leq j \leq 5)$ of the 
polyhedron $P$$)$. Then  $f_0(P) \leq 9$, and one of the following 
$4$ situations holds:

{\rm (i)} $f_0(P) = 8$, and one has the primitive relations
\[ v_0 + v_7 = v_1, \;\; v_0 + v_6 = 0; \]

{\rm (ii)} $f_0(P) = 8$, and one has the primitive relations
\[ v_0 + v_7 = v_3, \;\; v_0 + v_6 = 0; \] 

{\rm (iii)} $f_0(P) = 8$, and one has the primitive relations
\[ v_0 + v_7 = v_1, \;\; v_0 + v_6 = v_2; \]

{\rm (iv)} $f_0(P) = 9$, and one has the primitive relations
\[ v_0 + v_7 = v_1, \;\; v_0 + v_6 = v_2, \;\; v_0 + v_8 = 0. \]
\label{case.b}
\end{theo}

We consider separately all cases \ref{case.b}(i)-(iv).

\begin{abs}
\label{case1} 
{\rm  Consider the case \ref{case.b}(i).  Without loss of 
generality we can that there exists 
the primitive relation $v_1 + v_2 = v_0$ (othewise 
we could start with  $v_1$ instead of $v_0$). Thus, 
$D_0$ is a ${\P}^1$-bundle over ${\P}^2$: 
\[ {\P}_{{\P}^2} ( {\cal O} \oplus {\cal O}(a)). \]
By \ref{3-dim.face}, the $3$-polyhedron $P_0$ can not contain a 
$2$-dimensional face of the type $F(2)$. Therefore, $\mid a \mid \leq 2$, 
i.e., $P_0$ is a Fano $3$-polyhedron. 
By \ref{un-comb}, the combinatorial type of $P$ is uniquely determined. 
We denote it by $H$ (this type correspond to $P_{13}^8$ from 
\cite{grun}). 
Now it is easy to find the primitive relations defining $P$:   
$v_1  + v_2 = v_0$, $v_0 + v_7  = v_1$, $v_1 + v_6  = v_7$, $v_2 + v_7  = 0$, 
$v_0 + v_6 = 0$ and}
\end{abs}

\begin{center}
\begin{tabular}{|c|c|c|c|c|c|} \hline
\FP & $H_1$ & $H_2$ & $H_3$ & $H_4$ & $H_5$ \\ \hline
$v_3 + v_4 + v_5 = $ & $2v_1$ & $v_0 + v_1$ & $2v_0$ & $v_1$ & $v_0$ \\ \hline
\end{tabular}
\end{center}

\begin{center}
\begin{tabular}{|c|c|c|c|c|c|} \hline
\FP &  $H_6$ & $H_7$ & $H_8$ & $H_9$ & $H_{10}$ \\ \hline
$v_3 + v_4 + v_5 = $ & $v_0 + v_2$ & $2v_2$ & $0$ & $v_2$ & 
$v_2 +v_6$  \\ \hline
\end{tabular}
\end{center}
\medskip

\begin{abs} 
{\label{case2}} 
{\rm  Consider the case \ref{case.b}(ii). By \ref{num-eff} and 
\ref{3-dim.face}, ones obtains that either 
$D_{0}$ is a Fano $3$-fold, or $D_{0}$ is isomorphic to 
\[ {\P}_{{\P}^1} ( {\cal O} \oplus {\cal O}(1)\oplus {\cal O}(1)). \]
If $D_{0}$ is a Fano $3$-fold, then the combinatorial type of 
$P$ is uniquely determined by \ref{un-comb}. 
We denote this combinatorial type 
by $I$ (this type correspond to $P_{17}^8$ from 
\cite{grun}). The primitive relations are 
$v_0 + v_7  = v_3$, $v_3 + v_6  = v_7$, $v_0 + v_6  = 0$, and

\begin{center}
\begin{tabular}{|c|c|c|c|c|c|c|c|} \hline
\FP & $I_1$ & $I_2$ & $I_3$ & $I_4$ & $I_5$ & $I_6$ & $I_7$ \\ \hline
$v_1 + v_2 =$ & $v_3$ & $v_0$ & $v_0$ & $v_7$ & $v_4$ & $v_3$ 
& $0$ \\ \hline
$v_3 + v_4 + v_5 =$ & $2v_0$ & $2v_0$ & $v_0 + v_1$ & $2v_0$ & $2v_0$ & 
$v_0$ & $2v_0$ \\ \hline
$v_4 + v_5 + v_7 =$ & $v_0$ & $v_0$ & 
$v_1$ & $v_0$ & $v_0$ & $0$ & $v_0$ \\ \hline
\end{tabular}
\end{center}
\medskip

\begin{center}
\begin{tabular}{|c|c|c|c|c|c|c|c|c|} \hline
\FP & $I_8$ & $I_9$ & $I_{10}$ & $I_{11}$ & $I_{12}$ & $I_{13}$ & $I_{14}$ 
 & $I_{15}$ \\ \hline
$v_1 + v_2 =$ & $v_6$ & $v_7$ & $v_0$ & 0 & $v_6$ & $0$ & $v_4$ & $v_6$ 
\\ \hline
$v_3 + v_4 + v_5 =$ 
& $v_0+v_1$ & $v_0$ & $v_0$ & $v_0+v_1$ & $v_0$ & $v_0$ & $v_0$ & 
$2v_0$ \\ \hline
$v_4 + v_5 + v_7 =$ 
& $v_1$ & $0$ & $0$ & $v_1$ & $0$ & $0$ & $0$ & $v_0$ \\ \hline
\end{tabular}
\end{center}
\medskip
If $D_{0}$ is not a Fano $3$-fold, i.e., $D_0 \cong 
{\P}_{{\P}^1} ( {\cal O} \oplus {\cal O}(1)\oplus {\cal O}(1))$,  
then we determine uniquely the combinatorial type of $P$ using  \ref{un-comb}
and  denote this combinatorial type by $J$ (this type correspond to 
$P_{21}^8$ from 
\cite{grun}). The primitive relations defining 
$P$ are 
 $v_3 + v_6  = v_7$, $v_0 + v_1 + v_2 = v_4 + v_5$, $v_4 + v_5 + v_6 = 
v_1 + v_2$, $v_0 + v_7 = v_3$, $v_0 + v_6  = 0$ and}

\begin{center}
\begin{tabular}{|c|c|c|} \hline
\FP & $J_1$ & $J_2$   \\ \hline
$v_3 + v_4 + v_5 = $ & $0$ & $v_0$  \\ \hline
$v_4 + v_5 + v_7 = $ & $v_6$ & $0$ \\ \hline
$v_1 + v_2 + v_3 = $ & $v_6$ & $0$ \\ \hline
$v_1 + v_2 + v_7 = $ & $2v_6$ & $v_6$ \\ \hline
\end{tabular}
\end{center}
\end{abs}
\medskip

\begin{abs}
{\rm Consider the 
case \ref{case.b}(iii).  It this case the vertex $v_1$ satisfies 
the same conditions as the vertex $v_0$ in \ref{case.b}(i) and 
we do not obtain new Fano $4$-polyhedra. }
\end{abs} 

\begin{abs} 
{\rm Consider the case \ref{case.b}(iv). 
By \ref{3-dim.face}, $D_{0}$ must be a Fano $3$-fold,  
otherwise $P_0$ 
could have at most one double vertex (see \ref{commun}). 
On the other hand, all  $6$ vertices $v_0, v_1, v_2, 
v_6, v_7, v_8$ are contained 
in a $2$-dimensional linear subspace of ${\R}^4$. Thus, $P$ defines a 
locally trivial toric bundle over ${\P}^2$ whose fiber is 
a Del Pezzo surface with Picard rank $4$. If we denote the combinatorial 
type of $P$ by $K$, then the primitive relations defining $P$ are 
 $v_0 + v_7  = v_1$, $v_1 + v_8 = v_7$, $v_0 + v_8  = 0$, $v_8 + v_2 = v_6$, 
$v_7 + v_6  =  v_8$, 
$v_1 + v_6  = 0$, $v_0 + v_6  = v_2$, $v_1 + v_2  = v_0$, $v_7 + v_2 = 0$ 
and }
\label{case4}
\end{abs} 

\begin{center}
\begin{tabular}{|c|c|c|c|c|} \hline
\FP & $K_1$ & $K_2$ & $K_3$ & $K_4$  \\ \hline
$v_3 + v_4 + v_5 = $ & $2v_0$ & $v_0 + v_1$ & $v_0$ & $$0$$  \\ \hline
\end{tabular}
\end{center}
\medskip

Thus, we come to the following theorem. 

\begin{theo} There exist  exactly $31$ different Fano $4$-polyhedra 
$P$ such that $f_0(P) \geq 8$ and $P$ has a vertex of valence $5$. 
\end{theo}

\begin{coro} Among the combinatorial types 
$P_1^8, \ldots, P^8_{25}$, there exist exactly $3$ types $P_{13}^8$, 
$P_{17}^8$, and $P^8_{21}$ which admit a realization by Fano $4$-polyhedra.
\end{coro}

\noindent 
{\em Proof.} 
According to the table in \cite{grun}, for every combinatorial type 
from $P_1^8, \ldots, P^8_{12}$, there always exists a vertex $v_i$ 
of valence $4$. By \ref{commun}, $P_i$ contains at most $1$ double
vertex. This implies that the number of vertices of $P$ could be at 
most $7$. Contradiction to $f_0(P) =8$. For every combinatorial type 
from $P_{13}^8, \ldots, P^8_{25}$, there always exists a vertex  of 
valence $5$. Thus, the statement follows from \ref{case1}-\ref{case4}. 
\hfill $\Box$
\bigskip

\subsection{Fano $4$-polyhedra with $8$ vertices} 

\begin{prop} Let $V'$ $($resp. $V''$$)$ be a $3$-dimensional smooth 
projective toric variety defined by a  fan $\Sigma'$ $($resp. by 
$\Sigma''$$)$ 
having the combinatorial type $c_1$ $($resp. $c_2$$)$. Assume that both 
$V'$ and $V''$ have numerically effective anticanonical 
class and both $V'$ and $V''$ are not Fano $3$-folds.  
Let $P'$ $($resp. $P''$$)$ be the convex hull of generators of 
all $1$-dimensional 
cones in $\Sigma'$ $($resp. in $\Sigma''$$)$. 
If both polyhedra $P'$ and $p''$ have only  $2$-dimensional faces the of the 
types 
$F(1)$ and $F(4)$, then the fans $\Sigma'$ and $\Sigma''$  
are  defined by the following primitive relations among 
vertices $n_i$ $( 1 \leq i \leq 6):$ 

{\rm (i)} Type $c_1$ $:$ $n_6 + n_3 = n_1 + n_2$, $n_1 + n_2 + n_5 = n_6$, 
$n_1 + n_2 + n_4 = \lambda n_3$, where $\lambda = 0,1,2;$

{\rm (ii)} Type $c_2$ $:$ $n_1 + n_2 = n_3 + n_6$, $n_3 +  n_5 = 0$, 
$n_4 + n_6 = (\mu -1)  n_3$, where $\mu  =0,1,2.$
\label{not-fano}
\end{prop}  

\noindent 
{\em Proof. } The statement follows from the explicit 
description of $3$-dimensional smooth projective toric varieties 
with the Picard number $3$ (see \cite{bat5}). We remark also that 
two reflexive $4$-polyhedra $P'(\lambda)$ and $P''(\mu)$ corresponding to 
different values of $\lambda$ and $\mu$ are isomorphic if and 
only if $\lambda = \mu$. 
\hfill $\Box$

\begin{prop} Let $P$ be a Fano $4$-polyhedron such that 
among its vertices there exists the primitive relation  
$v_1 + v_2 = v_0$, where $D_{0}$ is a toric ${\P}^1$-bundle over 
a toric Del Pezzo surface with the Picard number $3$ or $4$. Then 
the toric variety $V'$ obtained from the contraction of the exceptional 
divisor $D_{0} \subset V(P)$ is again a toric Fano $4$-fold. 
\label{c-fano}
\end{prop}

\noindent 
{\em Proof. } Denote by $D_{[1,2]}$ the toric Del Pezzo 
surface on $V'$ obtained as the image of $D_{0}$. It is sufficient 
to prove that for any effective $1$-cycle $C$ on $D_{[1,2]}$ the intersection 
number $K_{V'} \cdot C$ is positive. Note that the cone of 
effective $1$-cycles on 
the Del Pezzo surface 
$D_{[1,2]}$ is always generated by exceptional rational curves of first 
kind. Thus, it is sufficient 
to assume that $C$ is a $1$-dimensional toric stratum on $D_{[1,2]}$ 
such that $C \cdot C = -1$. 

Let  $\R_{\geq 0} 
\langle v_1,v_2, v_3 \rangle \in \Sigma'$ be the $3$-dimensional 
cone of the fan $\Sigma'$
corresponding  to the 
stratum $C \subset V'$. Consider two $4$-dimensional cones 
$\R_{\geq 0} 
\langle v_1,v_2, v_3, v_4 \rangle$  and $\R_{\geq 0} 
\langle v_1,v_2, v_3, v_5\rangle$ in $\Sigma'$ having   
 $\R_{\geq 0} \langle v_1,v_2, v_3 \rangle$  as a common face. 
It  suffices to prove that in  the linear relation
\[ v_5 + v_4 + x_1 v_1 + x_2 v_2 + x_3 v_3 = 0 \]
the coefficients satisfy the inequiality 
$ x_1 + x_2 + x_3 \leq -1$. Since $C$ is an exceptional curve, 
we get $x_3 = -1$. Since the anticanonical divisor is ample on $V(P)$, 
we obtain that 
the coefficients of  two linear relations 
\[ v_5 + v_4 + y_1 v_0 + y_2 v_2 + y_3 v_3 = 0 \]
and 
\[ v_5 + v_4 + z_1 v_1 + z_2 v_0 + z_3 v_3 = 0, \]
satisfy the inequalities $ y_1 + y_2 + y_3 \leq -1$ and 
$z_1 + z_2 + z_3 \leq -1$. Combining these two relations with 
$v_1 + v_2 = v_0$,  we obtain 
\[ (y_1, y_2, y_3) = (x_1, x_2 - x_1, x_3)\;\;\; \mbox{\rm and } 
\;\;\; 
(z_1, z_2, z_3 ) = (x_1-x_2, x_2, x_3). \]
Thus, using $x_3 = -1$, we obtain  two inequalities  
$x_1 \leq 0$, $x_2 \leq 0$.
\hfill $\Box$

\begin{prop} Let $P$ be a Fano $4$-polyhedron with $8$ vertices 
such that there exists a primitive relation of the type $v_i + v_j = 0$. 
Then 
the following statements hold:

{\rm (i)} If $\Sigma_i(P)$ has the combinatorial type $c_1$ and 
$\Sigma_j(P)$ has the combinatorial type $c_2$, then $P$ has the combinatorial 
type $P_{26}^8$; 

{\rm (ii)} If both $\Sigma_i(P)$ and $\Sigma_j(P)$ have 
the combinatorial type $c_1$, then $P$ has the combinatorial 
type $P_{17}^8$; 

{\rm (iii)} If both $\Sigma_i(P)$ and $\Sigma_j(P)$ have 
the combinatorial type $c_2$, then $P$ has the combinatorial 
type $P_{34}^8$. 
\label{sym8}
\end{prop}

\noindent 
{\em Proof.} The statements follow from \ref{un-comb}. 
\hfill $\Box$
  
\begin{prop} Among the combinatorial types 
$P_{27}^8, \ldots , P_{33}^8$ only the type $P_{28}^8$ admits a 
realization by a Fano $4$-polyhedron.
\label{P8-28(1)}
\end{prop}

\noindent 
{\em Proof. } We shall use the same numeration of vertices as in 
\cite{grun}. 

Assume that $P$ has the combinatorial type $P_{27}^8$. Then 
$\{ v_6, v_8 \}$ is a primitive collection. The combinatorial type 
of $\Sigma_6(P)$ is $c_2$.  The combinatorial type 
of $\Sigma_8(P)$ is $c_1$. Since ${\cal V}(P)$ contains no vertex 
$v_i \neq v_6$ such that $\Sigma_i(P)$ has the combinatorial type 
$c_2$ or $d_2$, it follows from \ref{except} that $v_6 + v_8 =0$ is the 
single  possibility for the corresponding 
primitive relation. This contradicts \ref{sym8}(i). 

By similar method we obtain contradictions to  \ref{sym8}(ii) for the 
combinatorial types $P_{29}^8$, $P_{31}^8$, $P_{32}^8$, $P_{33}^8$,   
because for these $4$ combinatorial types there is no any vertex $v_i$ 
such that $\Sigma_i(P)$ has the combinatorial type $c_2$ or $d_2$. 

Finally, assume that $P$ has the combinatorial type $P^8_{30}$. 
Then ${\cal V}(P)$ contains only one primitive collection 
$\{ v_1, v_8 \}$ consisting of $2$ elements. We note that  
both $\Sigma_1(P)$, $\Sigma_8(P)$ have the combinatorial type 
$c_1$, there is no vertex $v_i \in {\cal V}(P)$ 
such that  $\Sigma_i(P)$ has the 
combinatorial type $c_2$,  and only $\Sigma_2(P)$ has the combinatorial 
type $d_2$. By \ref{sym8}(ii) and \ref{except}, the single 
possibility for the primitive relation is 
$v_1 + v_8 = v_2$. By \ref{c-fano}, the contraction of  the 
divisor $D_{2}$ yields  a toric Fano $4$-fold $V'$. The 
Fano $4$-polyhedron $P'$ corresponding to $V'$ has $7$ vertices, and any 
two vertices  $v_i, v_j \in P'$ are connected by  the edge  $[ v_i, v_j ]$ 
of $P'$. This contradicts \ref{poly7}. 
\hfill $\Box$

\begin{prop}
There exist exactly $2$ different Fano $4$-polyhedra having the 
combinatorial type $Z:= P_{28}^8$.  
\label{P8-28(2)}
\end{prop}

\noindent 
{\em Proof.} We shall use the same numeration of vertices of $P$ 
as in \cite{grun}. Then $\{v_1, v_8 \}$ and $\{v_5, v_7 \}$ 
are primitive collections. By \ref{sym8}, $\Delta(\{v_1, v_8 \}) = 
\Delta(\{v_5, v_7 \}) =1$.   By \ref{except},   
only $v_4$ and $v_6$ could be values for the sums $v_1 + v_8$ 
and $v_5 + v_7$. Since each edge $[v_1, v_6]$, $[v_5,v_4]$ is  
contained in exactly four $3$-dimensional 
faces of $P$ and  each edge 
$[v_1, v_4]$, $[v_5, v_6]$ is   contained in exactly 
five $3$-dimensional faces of $P$, the corresponding primitive 
relations are 
$v_1 + v_8 = v_4$, $v_5 + v_7 = v_6$. 
In order to find possibilities for other primitive relations, we 
apply \ref{c-fano} and obtain that the contraction $V'$ of $D_4$ or 
$D_6$ is again a Fano $4$-fold corresponding to 
a Fano $4$-polyhedron $P'$ with $7$ vertices having 
the combinatorial type $G$ (only $4$-polyehdra 
$G_4,G_6$) are possible). 
Using our classification 
in 3.1, we come to the following primitive relations:
$v_1 + v_2 + v_5 = 0$, $v_1 + v_2 + v_6 = v_7$, 
$v_2 + v_4 + v_5 = v_8$, $v_2 + v_4 + v_6 = v_7 + v_8$ and 

\begin{center}
\begin{tabular}{|c|c|c|} \hline
\FP & $Z_1$ & $Z_2$   \\ \hline
$v_3 + v_8 + v_7 = $ & $0$ & $v_2$  \\ \hline
$v_3 + v_4 + v_6 = $ & $v_1 + v_5$ & $0$ \\ \hline
$v_3 + v_4 + v_7 = $ & $v_1$ & $v_1 + v_2$ \\ \hline
$v_3 + v_6 + v_8 = $ & $v_5$ & $v_2 + v_5$ \\ \hline
\end{tabular}
\end{center}
\hfill $\Box$

\begin{prop} 
\label{P8-34}
There exist exactly $13$ different Fano 
$4$-polyhedra having the combinatorial type $L := P_{34}^8$. 
The primitive relations are  $v_1 + v_8  = 0$ and 

\begin{center}
\begin{tabular}{|c|c|c|c|c|c|c|c|} \hline
\FP & $L_1$ & $L_2$ & $L_3$ & $L_4$ & $L_5$ & $L_6$ & $L_7$ \\ \hline
$v_2 + v_3 =$ & $v_1$ & $v_1$ & $v_1$ & $v_1$ & $0$ & $0$ 
& $v_1$ \\ \hline
$v_4 + v_5 =$ & $v_1$ & $v_3$ & $v_1$ 
& $v_3$ & $v_3$ & $v_3$ & 
$0$ \\ \hline
$v_6 + v_7 = $ & $v_1$ & $v_3$ & 
$v_4$ & $v_4$ & $v_3$ & $v_4$ & $v_4$ \\ \hline
\end{tabular}
\end{center}
\medskip

\begin{center}
\begin{tabular}{|c|c|c|c|c|c|c|} \hline
\FP & $L_8$ & $L_9$ & $L_{10}$ & $L_{11}$ & $L_{12}$ & $L_{13}$ 
 \\ \hline
$v_2 + v_3 =$ & $0$ & $0$ & $v_1$ & $0$ & $v_1$ & $v_1$  \\ \hline
$v_4 + v_5 =$ & $0$ & $0$ & $v_3$ & $v_3$ & $v_8$ & $v_1$  \\ \hline
$v_6 + v_7 = $ & $0$ & $v_4$ & $v_2$ & $v_2$ & $v_4$ & $v_8$ \\ \hline
\end{tabular}
\end{center}
\medskip
\end{prop}

\noindent 
{\em Proof.} The statement follows directly from the results in 
\cite{bat5}, because  the combinatorial structure of $P$ 
is described by  pairwise disjoint primitive collections.  
 
\hfill $\Box$

\begin{prop} 
There exist exactly $4$ different Fano $4$-polyhedra having  
the combinatorial type type $M:= P_{26}^8$ and a pair of 
centrally symmetric vertices $v_1, v_8$. 
\label{P8-26(1)}
\end{prop}

\noindent 
{\em Proof. } By \ref{z-sym} and \ref{un-comb}, if $v_1$ and $v_8$  
are centrally symmetric 
vertices, then $P_1 := \pi_1(P)$ is not a Fano polyhedron. Moreover, 
$P_1$ contains exactly one $2$-dimensional face of the type $F(4)$. 
Therefore all possibilities for $P_1$ are described in \ref{not-fano}. 
Using this descriptions of $P_1$, we come to the following 
primitive relations: 
 $v_1 + v_8  = 0$, $v_1 + v_2 + v_3 = v_4 + v_6$, $v_4 + v_6 + v_8
 =v_1 + v_2$   and

\begin{center}
\begin{tabular}{|c|c|c|c|c|} \hline
\FP & $M_1$ & $M_2$ & $M_3$ & $M_4$   \\ \hline
$v_4 + v_5 =$ & $0$ & $v_1$ & $v_1$ & $v_1$   \\ \hline
$v_6 + v_7 = $ & $0$ & $v_1$ & $v_5$ & $0$   \\ \hline
$v_2 + v_3 + v_5 = $ & $v_6 + v_8$ & $v_6$ & $v_6$ & $v_6$ \\ \hline
$v_2 + v_3 + v_7 = $ & $v_4 + v_8$ & $v_4$ & $0$ & $v_4 + v_8$  \\ \hline
\end{tabular}
\end{center}

\hfill $\Box$

\begin{prop} 
There exists exactly one Fano $4$-polyhedron having  
the combinatorial type $M := P_{26}^8$ and containing no pair of 
centrally symmetric vertices. 
\label{P8-26(2)}
\end{prop}

\noindent 
{\em Proof.} 
Using the numeration of vertices of $P$ as in \cite{grun}, we obtain the 
following primitive collections
\[ \{v_1, v_8\},\; \{v_4, v_5\},\; \{v_6, v_7\},\; \{v_1, v_2, v_3\},\;
\{v_2, v_3, v_5\},\;\{v_2, v_3, v_7\},\; \{v_4, v_6, v_8\}. \]
The fans  $\Sigma_4(P)$, $\Sigma_6(P)$, 
$\Sigma_8(P)$ have the 
combinatorial type $c_1$. The fans $\Sigma_1(P)$, $\Sigma_5(P)$, 
$\Sigma_7(P)$ have the combinatorial type $c_2$. The fans 
$\Sigma_2(P)$ and  
$\Sigma_3(P)$ have the combinatorial type $d_5$ (see the notations 
in  \cite{grun}). 

By \ref{except}, the sum $v_1 + v_8$ must be equal to either  
$v_5$ or $v_7$. Without loss of generality, we assume that $v_1 + v_8 = v_5$. 
Applying again \ref{except}, we see that the sum $v_4 + v_5$ must  
be equal to either $v_1$ or $v_7$. If $v_4 + v_5 = v_1$, 
then $v_4 + v_8 = 0$, contradiction to our assumtion on $P$. 
Thus, the single  possibility is 
$v_4 + v_5 = v_7$. By similar method, we obtain $v_6 + v_7 = v_1$. 
This implies that  $v_4 + v_6 + v_8 = 0$. Using already known primitive 
relations, we can restrict  possibilies for the fan $\Sigma_2(P)$ having 
the type $d_5$ such that the convex hull 
of the generators of  all $1$-dimensional cones 
in $\Sigma_2(P)$ is a reflexive  polyhedron 
whose  $2$-dimensional faces isomorphic  
to $F(1)$ or   $F(4)$. As a result, we 
obtain the required  primitive relations which are uniquely defined up to 
the cyclic permutation of the primitive pairs of vertices 
\[ \{v_1, v_8\} \rightarrow  \{v_4, v_5\} \to \{v_6, v_7\} \to  
\{v_1, v_8\}. \]
So we obtain the following primitive relations defining the polyhedron 
$M_5$: 
$v_1 + v_8 = v_5$,  $v_4 + v_5  = v_7$, $v_6 + v_7  = v_1$, 
$v_1 + v_2 + v_3 = v_6$, $v_2 + v_3 + v_5 = v_6 + v_8$,       
$v_2 + v_3 + v_7 = 0$, $v_4 + v_6 + v_8 = 0$. 

\hfill $\Box$

\begin{prop} 
The combinatorial types $P_{35}^8$, $P_{36}^8$, $P_{37}^8$ do not 
admit a realization by a Fano $4$-polyhedron $P$. 
\label{567}
\end{prop}

\noindent
{\em Proof.} For  combinatorial types $P_{35}^8$, $P_{36}^8$, $P_{37}^8$  
all vertices have valence $7$. Therefore, 
${\cal V}(P)$ does not contain a primitive collection consisting  
of exactly two vertices. This contradicts \ref{two.vert}
\hfill $\Box$

Collecting together all our results, we come to the following:

\begin{theo} There exist exactly $47$ different Fano $4$-polyhedra 
with $8$ vertices. Among $37$ all possible combinatorial types of simplicial 
convex $4$-polyhedra with $8$ vertices only $6$ ones $P_{13}^8$,  
$P_{17}^8$, $P_{21}^8$, $P_{26}^8$, $P_{28}^8$, $P_{34}^8$ 
admit a realization by a Fano $4$-polyhedron. 
\label{8-vert}
\end{theo}

\subsection{Fano $4$-polyhedra having a vertex of valence 6} 

\begin{prop} 
Let $P$ be a Fano $4$-polyhedron with at least 
$9$ vertices which contains  a vertex $v_i$ of valence 
$6$ and does not contain vertices of 
valence $\leq 5$. Then either $P$ contains  a vertex $v_j$ such that the fan 
$\Sigma_j(P)$ has 
the combinatorial type $c_2$, or 
$f_0(P)=9$ and $P$ is isomorphic to the following 
Fano $4$-polyhedron defined by the primitive relations: 
$v_0 + v_7 = 0$, $v_0 + v_8 = v_1$, $v_3 + v_5 = v_4$, $v_4 + v_6 = v_5$, 
$v_1 + v_7 = v_8$, $v_3 + v_6 = 0$, $v_1 + v_2 + v_5 = 
v_0 + v_6$, $v_1 + v_2 + v_4 = v_0$, $v_2 + v_5 + v_8 = v_6$, 
$v_1 + v_2 + v_8 = 0$. 
\label{single}
\end{prop} 

\noindent
{\em Proof.} There exist exactly two different possibilities for 
the combinatorial type of $\Sigma_i(P)$: $c_1$ and $c_2$. 
In the second case we can put $j = i$. Now assume that $\Sigma_i(P)$
has the type $c_1$. Putting $i =0$ and using the numeration
for the vertices joined with $v_0$ as in \ref{stereo}, we obtain at least 
two primitive collections  
$\{ v_0, v_7\}$ and $\{v_0, v_8 \}$. So, at least one 
of two points $\pi_0(v_7), \pi_0(v_8) \in P_0 = \pi_0(P)$ is 
a double vertex of $P_0$. Assume that $\pi_0(v_8)$ is a double vertex. 
Then, by \ref{double-2},  
we have a primitive relation $v_0 + v_8 = v_k$ for some 
vertex $v_k \in P$ $(1 \leq k \leq 6)$. If the corresponding vertex 
$n_k$ in the graph $c_1$ has valence $3$ (i.e., $k \in \{3,6\}$), then $v_k$ has valence $5$ 
(contradiction). If $n_k$ in the graph $c_1$ has valence $4$ 
(i.e., $k \in \{4,5\}$), then 
$v_k$ has valence $6$. By \ref{ex-ray2}(iii), $\Sigma_k(P)$ has the 
combinatorial type $c_2$. The single remaining possibility 
is $k \in \{1, 2 \}$. 

Assume that $k =1$. By \ref{commun}, $\pi_0(v_1)$ is the unique 
double vertex of $P_0$. Therefore, we obtain that $f_0(P) \leq 9$, i.e., 
$f_0(P) = 9$. Moreover, by \ref{c-fano}, the contraction of 
$D_1$ is again a Fano $4$-fold $V'$ corresponding to a Fano $4$-polyhedron 
$P'$ having $8$ vertices. Using our results in 3.2 and 3.3, we obtain 
that the combinatorial type of $P'$ equals $I$ and $P'$ must be isomorphic 
to $I_{13}$. The last  fact implies the above primitive relations 
describing  $P$. 
\hfill $\Box$

\begin{prop} 
There exist exactly $20$ different Fano $4$-polyhedra $P$ having 
$9$ vertices and containing a vertex $v_i$ such that the fan 
$\Sigma_i(P)$ has 
the combinatorial type $c_2$. 
\label{c2-class}
\end{prop}

\noindent 
{\em Proof.} We assume that $i = 0$, i.e.,  $v_0 \in  P$ is 
the vertex of valence $6$. We choose  a numeration of 
vertices $v_j$ ($1 \leq j \leq 6$) of $P$ 
in such a way that it corresponds to the numeration of
 the vertices $n_j$ ($1 \leq j \leq 6$) in the  graph $c_2$. 
Then  $\{ v_0, v_7 \}$ and $\{ v_0, v_8 \}$ are 
primitive collections. There exist  the following three cases 
for the corresponding 
primitive relations:
\medskip

\noindent
{\sc  Case I}: $P_0$ is a Fano $3$-polyhedron and the  
primitive relations are 
\[ v_0 + v_8 = 0,\;\; v_0 + v_7 = v_k\; (1 \leq k \leq 6).\]
Without loss of generality, we can assume $k =1$. 
Then the combinatorial type of $P$ is determined by 
\ref{un-comb} and $P$ is determined by the primitive relations 
 $v_1 + v_2  = v_0$, $v_1 + v_8  = v_7$, 
$v_2 + v_7  = 0$ and

\begin{center}
\begin{tabular}{|c|c|c|c|c|c|c|} \hline
\FP & $Q_1$ & $Q_2$ & $Q_3$ & $Q_4$ & $Q_5$ & $Q_6$ \\ \hline
$v_3 + v_5 =$ & $v_1$ & $v_1$ & $v_0$ & $v_1$ & $v_0$ & $$0$$  \\ \hline
$v_4 + v_6 =$ & $v_1$ & $v_3$ & $v_0$ & $v_0$ & $v_3$ & $v_1$  \\ \hline
\end{tabular}
\end{center}
\medskip

\begin{center}
\begin{tabular}{|c|c|c|c|c|c|c|} \hline
\FP & $Q_7$ & $Q_8$ & $Q_9$ & $Q_{10}$ & $Q_{11}$ & $Q_{12}$ \\ \hline
$v_3 + v_5 =$ & $v_0$ & $$0$$ & $v_0$ & $$0$$ & $$0$$ & $v_1$  \\ \hline
$v_4 + v_6 =$ & $v_7$ & $v_0$ & $v_2$ & $v_3$ & $$0$$ & $v_2$  \\ \hline
\end{tabular}
\end{center}
\medskip

\begin{center}
\begin{tabular}{|c|c|c|c|c|c|} \hline
\FP & $Q_{13}$ & $Q_{14}$ & $Q_{15}$ & $Q_{16}$ & $Q_{17}$  \\ \hline
$v_3 + v_5 =$ & $v_2$ & $v_2$ & $$0$$ & $v_0$ & $v_2$   \\ \hline
$v_4 + v_6 =$ & $v_2$ & $v_3$ & $v_2$ & $v_8$ & $v_8$   \\ \hline
\end{tabular}
\end{center}
\medskip

\noindent
{\sc Case II:}  $P_0$ is not a Fano $3$-polyhedron and the  
primitive relations are 
\[ v_0 + v_8 = 0,\;\; v_0 + v_7 = v_k\; (1 \leq k \leq 6).\]
Then $P_0$ is isomorphic to one of $3$ polyhedra from \ref{not-fano} and 
the combinatorial type of $P$ is uniquely determined by 
\ref{un-comb}. By \ref{commun}, only $\pi_0(v_4)$ or $\pi_0(v_5)$ 
can be double a vertex of $P_0$. Assume that $v_4$ is a double vertex, i.e., 
$k =4$. Then we obtain the primitive relations  
 $v_0 + v_1 + v_2 = v_3 + v_6$, $v_3 + v_6 + v_8 = v_1 + v_2$,  
$v_0 + v_7  = v_4$, $v_4 + v_8 = v_7$, $v_0 + v_8 =0$.  
By \ref{un-comb}, $\{v_7,v_6\}$ and $\{v_4,v_6\}$ are primitive collections. 
Therefore,  $\pi_0(v_4)$ and $\pi_0(v_6)$ must be centrally symmetric. 
So, the lattice points  $v_0, v_4, v_6, v_7, v_8$ are vertices of a 
$2$-dimensional Fano polyhedron and one of two pairs 
$\{ v_4, v_6 \}$ or $\{ v_7, v_6 \}$ consists of centrally symmetric 
vertices. Without loss of generality we shall assume that 
$ v_7 + v_6 = 0$. Then, using the classification 
of Fano $4$-polyhedra of the combinatorial type $P_{26}^8$, we obtain 
the primitive relations:  $v_1 + v_2 + v_7 = v_3 + v_8$, 
$v_1 + v_2 + v_4 = v_3$ and  

\begin{center}
\begin{tabular}{|c|c|c|c|} \hline
\FP & $R_1$ & $R_2$ & $R_3$  \\ \hline
$v_3 + v_5 = $ & $v_4$ & $v_0$ & $0$  \\ \hline
$v_1 + v_2 + v_5 = $ & $0$ & $v_6$ & $v_6+v_8$  \\ \hline
\end{tabular}
\end{center}
\medskip

\noindent
{\sc Case III:}  The primitive relations are 
\[ v_0 + v_8 = v_l,\;\; v_0 + v_7 = v_k\; (1 \leq l \neq k \leq 6).\]
By \ref{commun}(ii), the pair $\{v_k, v_l\}$ must coincide 
with one of the pairs $\{ v_1, v_2 \}$, $\{ v_3, v_5 \}$, $\{ v_4, v_6 \}$. 
Assume that $k =1$ and $l =2$. We note that $P_0$ can contain only faces 
of the type $F(1)$ or $F(4)$ (otherwise there would exist a primitive 
collection of the type $v_0  + v_{i_1} + v_{i_2} = 2v_j$ or 
$v_0  + v_{i_1} + v_{i_2} + v_{i_3} = 3v_j$, and this would imply, by 
\ref{ex-ray2}, that $v_i$ is a vertex of valence $\leq 5$). 
Furthermore, by \ref{not-fano}, we obtain that $P_0$ must 
a Fano $3$-polyhedron. 

If $\pi_0(v_1) + \pi_0(v_2) = \pi_0(v_i) \neq 0$ for 
some vertex $v_i \in {\cal V}(P)$, then we obtain two primitive 
relations $v_1 + v_8 = v_i$ and $v_2 + v_7 = v_i$. On the other hand, by 
\ref{except}, $\Sigma_1(P)$ and $\Sigma_2(P)$ have the combinatorial type 
$c_2$. This implies that $\{v_1, v_2\}$ is also a primitive collection. 
Considering the  $\pi_0$-projection, we obtain the unique
possibility: $v_1 + v_2 = v_i$. But the last equation contradicts
$v_1 + v_8 = v_i$ and $v_2 + v_7 = v_i$.

It remains the single possibility $\pi_0(v_1) + \pi_0(v_2) = 0$.  
Then all $5$ vertices $v_1$, 
$v_2$, $v_0$, $v_7$, $v_8$ are contained in a $2$-dimensional linear 
subspace in ${\R}^4$ and form a $2$-dimensional Fano polyhedron. 
Thus, the vertex $v_1$ satisfies the 
same conditions as the vertex $v_0$ in Case I.
We have shown that the considered  case reduces to the case I, 
so we obtain no new 
Fano $4$-polyhedra.  
\hfill $\Box$

\begin{prop} 
There exist exactly $8$ different Fano $4$-polyhedra $P$ having 
at least $10$ vertices and containing a vertex $v_i$ such that 
$\Sigma_i(P)$ has the combinatorial type $c_2$. 
\label{c2-class2}
\end{prop}

\noindent 
{\em Proof. }
The primitive relations are 
$v_1 + v_3  = v_2$, $v_2 + v_4 = v_3$, $v_1 + v_4 = 0$,  $v_3 + v_5  = v_4$, 
$v_4 + v_6  = v_5$, $v_2 + v_5 = 0$,  $v_1 + v_5  = v_6$, $v_2 + v_6 = v_1$, 
$v_3 + v_6  = 0$

\begin{center}
\begin{tabular}{|c|c|c|c|c|c|c|c|c|} \hline
\FP &  $U_1$ & $U_2$ & $U_3$ & $U_4$ & $U_5$ & $U_6$ & $U_7$ & $U_8$ \\ \hline
$v_8 + v_7 = $ & $v_1$ & $v_1$ & $v_1 $ & $0$ &  $0$ & $v_1$ & $v_1$ & $v_1$ 
\\ \hline
$v_9 + v_{10} = $ & $v_1$ & $v_8$ & $v_2$ & $v_8$ & $0$ & $0$ & $v_3$ & $v_4$ 
\\ \hline
\end{tabular}
\end{center}
\hfill $\Box$
\bigskip

\subsection{Fano $4$-polyhedra whose all vertices have valence 
$\geq 7$}

\begin{prop}  
For a Fano $4$-polyhedron $P$, one has 
\[ 4f_0(P) \geq f_1(P). \]
\label{inequal-1}
\end{prop}
 
\noindent 
{\em Proof.} The statement follows from the general inequality in 
\ref{weights-1} together with  the Dehn-Sommerville relation 
$f_2(P) = 2(f_1(P) - f_0(P))$. 
\hfill $\Box$

\begin{coro}
Assume that any two vertices of a Fano $4$-polyhedron $P$ are joined by a 
$1$-dimensional face, then $f_0(P) \leq 9$. 
\label{leq9}
\end{coro}

\noindent 
{\em Proof.}
Note that 
\[ f_1(P) = \frac{f_0 (P) ( f_0(P) -1 )}{2}.\]
By \ref{inequal-1}, one has $f_0(P) \leq 9$.
\hfill $\Box$

\begin{theo}
Let $P$ be a Fano $4$-polyhedron 
such that all vertices of $P$ have valence $\geq 7$. Then there 
always exists a primitive collection containing exactly  two vertices of 
$P$.
\label{two.vert}
\end{theo}

\begin{lem}
Assume that any two vertices of a Fano $4$-polyhedron $P$ are joined by a 
$1$-dimensional face, then all primitive 
relations corresponding to primitive collections of degree $1$ have the 
form  
\[ v_{i_1} + v_{i_2} + v_{i_3} = v_{j_1} + v_{j_1}. \]
\end{lem}

\noindent 
{\em Proof.} 
Assume that there exists a primitive relation of 
the type  
\[ v_{i_1} + v_{i_2} + v_{i_3} + v_{i_4} = 3v_{j}. \]
Then, by \ref{ex-ray2}(iii),  $D_{j}$ is isomorphic to ${\P}^3$, i.e., 
$v_j$ has valence $4$. Contradiction. 

Assume that there exists a primitive relation of 
the type  
\[ v_{i_1} + v_{i_2} + v_{i_3} = 2v_{j}, \]
Then, by \ref{ex-ray2}(iii),  $D_{j}$ is a toric ${\P}^2$-bundle over 
${\P}^1$, i.e., $v_j$ has valence $5$. Contradiction. 

\hfill $\Box$

\begin{lem}
Assume that any two vertices of a Fano $4$-polyhedron are joined by a 
$1$-dimensional face and 
\[ v_{1} + v_{2} + v_{3} = v_{4} + v_{5} \]
is a primitive relation. Then for any other   primitive 
relations of the type 
\[ v_{i_1} + v_{i_2} + v_{i_3} = v_{j_1} + v_{j_2}, \]
one has 
\[ \{ v_{j_1}, v_{j_2} \} 
\not\subset \{ v_{1}, v_{2}, v_{3},v_{4},v_{5} \}. \]
\label{two.one}
\end{lem}

\noindent
{\em Proof. } Assume that $\{ v_{j_1}, v_{j_2} \} 
\subset \{ v_{1}, v_{2}, v_{3},v_{4},v_{5} \}$. Then there exists 
a vertex $v_i \in \{ v_{1}, v_{2}, v_{3} \}$ which is not 
contained in $\{ v_{j_1}, v_{j_2} \}$. Without loss of generality 
we assume that $v_i = v_3$.   
Then, by \ref{ex-ray2}(i), the edge $[v_4, v_5] \subset P$ 
is contained in the $3$-dimensional face  $[v_1,v_2,v_4,v_5]$.  
On the other hand,  by \ref{ex-ray2}(i), the edge $[ v_{j_1}, v_{j_2}]$  
is contained in exactly three $3$-dimensional faces of $P$: 
$$[v_{i_1},v_{i_2}, v_{j_1}, v_{j_2}], [v_{i_1},v_{i_3}, v_{j_1}, v_{j_2}], 
[v_{i_2},v_{i_3}, v_{j_1}, v_{j_2}].$$ 
This shows that  $[v_1,v_2,v_4,v_5]$ must coincide with one 
of these $3$ faces. Assume that    
$$ [ v_{i_1},v_{i_2}, v_{j_1}, v_{j_2} ] = 
[ v_1,v_2,v_4,v_5].  $$  
Then our two primitive relations imply $v_3 = v_{i_3}$. 
Contradiction. 

\hfill $\Box$

\begin{lem}
Assume that any two vertices of a Fano $4$-polyhedron are joined by a 
$1$-dimensional edge. Then the number 
of different primitive relations  of depth $1$ is not more 
than  $ f_3(P)/3$. 
\label{prim.num}
\end{lem}

\noindent 
{\em Proof.} Let $v_1 , v_2 , v_3 , v_4 $ be vertices of 
a $3$-dimensional face of $P$. By \ref{two.one}, there exists at most 
one primitive relation of the type 
\[ v_{i_1} + v_{i_2} + v_{i_3} = v_{j_1} + v_{j_2} \]
such that 
\[ \{ v_{j_1},  v_{j_2} \} \subset \{ v_1 , v_2 , v_3 , v_4 \}. \]
On the other hand, for any primitive relation of depth $1$ as above, 
 there exist exactly three $3$-dimensional faces of $P$ containing 
 $[ v_{j_1},  v_{j_2} ]$. 
\hfill $\Box$

\begin{lem} 
Assume that any two vertices of a Fano $4$-polyhedron are joined by a 
$1$-dimensional face.  Then the number 
of $2$-dimensional faces of $P$ having weight $-1$ is not more 
than  $18$  $($resp. $27$$)$  if $f_0(P) = 8$ $($resp. if $f_0(P) = 9$$)$. 
\label{faces.2}
\end{lem}

\noindent 
{\em Proof.} If $f_0(P) = 8$, then $f_3(P) = 20$. Thus, according 
to \ref{prim.num} the number of primitive relations of depth $1$ is not 
more than $[20/3] =6$. On the other hand, 
by \ref{ex-ray2} and by arguments in the proof of 
\ref{weights-1}, every $2$-dimensional face 
of $P$ having weight $-1$ appears 
from a primitive relation of degree  $1$. 
Moreover, every primitive relation of the type 
\[ v_{i_1} + v_{i_2} + v_{i_3} = v_{j_1} + v_{j_2} \] 
gives rise to exactly three $2$-dimensional faces of $P$ 
\[ [ v_{i_1}, v_{j_1}, v_{j_2} ],\; 
[ v_{i_2}, v_{j_1}, v_{j_2} ],\; 
 [ v_{i_3}, v_{j_1}, v_{j_2} ] \]
having the weight $-1$. Thus, the number of these $2$-dimensional faces 
$\leq 18$. 

If $f_0(P) = 9$, then $f_3(P) = 27$. By similar arguments, we obtain 
the number of $2$-dimensional faces of $P$ having weight $-1$ is not more 
than  $27$.   

\hfill $\Box$

\medskip 
\noindent
{\em Proof of theorem \ref{two.vert}} 
Assume that any two vertices of 
a Fano polyhedron $P$ are joined by a $1$-dimensional face. 
Using Dehn-Sommerville equalities and arguments in  
the proof of \ref{weights-1}, we obtain  that the total weight of 
the Fano polyhedron $P$ 
equals 
\[ w(P) = 18 f_0(P) - 6f_1(P). \]
Denote by $f_2^*(P)$ the number of $2$-dimensional 
faces of $P$ having the weight $-1$. Since the weight of 
every $2$-dimensional face of $P$ is at least $-1$, we obtain  
$$w(P) + f_2^{*}(P) \geq 0.$$ 

If $f_0(P) = 8$, then $w(P) =  18 \cdot 8 - 6 \cdot 28 = -24$. By 
\ref{faces.2}, $w(P) + f_2^{*}(P) \leq -6$. Contradiction. 

If $f_0(P) = 9$, then $w(P) =  18 \cdot 9 - 6 \cdot 36 = -54$. By 
\ref{faces.2}, $w(P) + f_2^{*}(P) \leq -27$. Contradiction. 
\hfill $\Box$

\begin{prop} Let $P$ be a Fano $4$-polyhedron. Assume  that 
for any primitive collection $\{ v_i, v_j \} \subset {\cal V}(P)$ 
the corresponding primitive relation has form $v_i + v_j = 0$. 
Then the following statements hold:  

{\rm (i)} $f_0 (P) \leq 10$;

{\rm (ii)} if $f_0(P) = 10$, then $P$ is isomorphic to the convex hull of 
the vectors   
\[ \pm e_1 , \; \pm e_2 , \; \pm e_3 , \; 
\pm e_4 , \; \pm(e_1  + e_2  + e_3  + e_4 ). \]

{\rm (iii)} if $f_0(P) = 9$, then $P$ is isomorphic to the convex hull of 
the vectors 
\[ \pm e_1 , \; \pm e_2 , \; \pm e_3 , \; \pm e_4 ,\;  
(e_1  + e_2  + e_3  + e_4 ) . \]
\label{symm-10}
\end{prop}

\noindent
{\em Proof. } 
(i) Since every vertex $v_i \in {\cal V}(P)$ 
is joined by a $1$-dimensional face 
with at least $f_0(P) -2$ other vertices of $P$, one has 
\[ 2f_1(P) \geq f_0(P) ( f_0(P) - 2 ). \]
On the other hand, $2f_1(P) \leq 8f_0(P)$ (see  \ref{inequal-1}).  
Thus, $f_0 \leq 10$. 

(ii) By \ref{two.vert},  there exists at least one primitive collection 
consisting of two vertices, e.g.,  $v_9$ and $v_{10}$. We shall assume 
that the corresponing primitive relation is  
$v_9 + v_{10} = 0$. Let $P_{10}$ be the $3$-dimensional projection of $P$ into 
${\R}^4 / {\R}<v_{10}>$.  By \ref{un-comb}, 
the combinatorial type of 
$\Sigma_{10}(P)$ and the polyhedron $P_{10}$ allow us to recover the 
combinatorial structure of the Fano $4$-polyhedron $P$. Since 
$\Sigma_{10}(P)$  consists of $12$ cones of dimension $4$, 
the polyhedron $P_{10}$ has at most six  $2$-dimensional faces of the 
type $F(4)$. On the other hand, there are $10$ pairs 
$\{ v_i', v_j' \}$ of vertices of $P_{10}$ such that the cone  
$\R_{\geq 0} v_i' + \R_{\geq 0} v_j'$ does not belong to  
the fan $\Sigma_{10}(P)$. By \ref{un-comb}, any such a pair 
$\{ v_i', v_j' \}$ gives rise to a primitive collection 
$\{ v_i, v_j , v_{10} \} \subset {\cal V}(P)$ 
unless $ v_i', v_j' $ are opposite vertices of some $2$-dimensional face of 
$P_{10}$ having the combinatorial type $F(4)$. Thus, among vertices 
$\{ v_1, \ldots, v_8 \}$ there exists at least $10 - 6 = 4$ primitive 
collections containing only 
two vertices. This shows that vertices  $\{ v_1, \ldots, v_8 \}$ can be 
divided into $4$ pairs of centrally symmetric vectors. Thus, $P$ is 
centrally symmetric Fano $4$-polyhedron having $10$ vertices. By 
results of Ewald and Voskresenski\^i \& Klyachko 
\cite{ewald,vosk.klyachko}, there exists a 
unique possibility for $P$. 

(iii) By \ref{two.vert},  there exists at least one primitive collection 
containing exactly two vertices, e.g.,  $v_8$ and $v_{9}$. Therefore, we have 
$v_8 + v_{9} = 0$. Let $P_9$ be the $3$-dimensional projection of $P$ into 
${\R}^4 / {\R}<v_{9}>$. Let us prove by the same method as 
in (ii) that among $\{ v_1, \ldots, v_7 \}$ there exist $3$ pairs 
of centrally symmetric vertices. By our assumption about $P$, it suffices to 
prove that there exist at least  $3$ primitive collections in 
$\{ v_1, \ldots, v_7 \}$ containing exactly $2$ vertices. 

Denote by $f(P_9)$ the number of $2$-dimensional faces of $P$ having 
the type $F(4)$. 
Since $\Sigma_9(P)$ 
contains exactly  $10$ cones of dimension $3$, 
we obtain $f(P_9) \leq 5$.  
On the other hand, there exist  
$6$ pairs $\{ v_i', v_j' \}$ of vertices of $P_9$ such that 
$\R_{\geq 0} v_i' + \R_{\geq 0} v_j'$ does not belong to  
the fan $\Sigma_{9}(P)$. By \ref{un-comb}, any such a pair 
$\{ v_i', v_j' \}$ gives rise to a primitive collection $\{ v_i, v_j, v_9 \}$ 
in ${\cal V}(P)$ 
unless $ v_i', v_j' $ are opposite vertices of some $2$-dimensional face of 
$P_9$ having the combinatorial type $F(4)$. Thus, among vertices 
$\{ v_1, \ldots, v_7 \}$ there exist at least $6 - f(P_9) \geq 1$ primitive 
collections.

Let us prove that $f(P_9) \neq 5,4$. 

Assume that $f(P_9) = 5$. Then $f(P_9) = f_2(P_9)$, and 
$f_1(P_9) = 10$. On the other 
hand,  every vertex $v_i' \in P_9$ is contained in at least three 
$1$-dimensional faces of $P_9$. Therefore, $f_1(P_9) \geq 7\cdot 3/2$. 
Contradiction. 

Assume that $f(P_9) = 4$. Then $f_2(P_9) = 6$, and $f_1(P_9) = 11$. For every 
vertex $v_i \in \{ v_1, \ldots, v_7 \} $, 
we denote by $r_i$ the number of $1$-dimensional 
faces of $P_9$ containing $\pi_9(v_i)$, and denote by $c_i$ the number of 
$2$-dimensional faces of $P_9$ having the type $F(4)$ and containing 
$\pi_9(v_i)$. 
By \ref{un-comb}, $6 - c_i - r_i$ is the number of primitive collections  
in ${\cal V}(P)$ of the type $\{ v_i, v_j \}$ where $v_j \in 
\{ v_1, \ldots, v_7 \} \setminus v_i$. By our assumption on $P$, 
$c_i + r_i$ must be $5$ or $6$. On the other hand, we have 
\[ \sum_{i =1}^7 r_i = 2 f_1(P_9) = 22, \;\; 
\sum_{i =1}^7 c_i = 4 f(P_9) = 16. \]
There exists a unique representation of $38$ as a sum of $7$ numbers 
which are equal to $5$ or 
$6$: 
\[ 38 = 5 + 5 + 5 + 5 + 6 + 6 + 6 =  \sum_{i =1}^7 (r_i +c_i). \] 
This representation determines  
the combinatorial type of $P_9$ uniquely. In particular, there always 
exists a vertex  $\pi_9(v_7) \in P_9$   
which is a common vertex of two triangles at 
the boundary of $P_9$. Since $f(P_9) =4$, there exist two pairs of 
centrally symmetric vertices among $\{ v_1, \ldots,v_6 \}$. We assume, 
for instance, that $v_1 + v_2 = v_3 + v_4 = 0$. Then, applying 
an authomorphism of $3$-dimensional integral lattice, we can transform 
$v_1', v_2', v_3' , v_4', v_7'$  to the form 
\[ v_1' = (1,0,0) = - v_2', \; \; v_3' = (0,1,0) = - v_4', \;\;
 v_7' = (0,0,1). \]
For  vertices $v_6'$ and $v_5'$ we have only the following $4$ possibilities 
\[ (1,1,-1),\; (-1,1,-1),\;(1,-1,-1),\;(-1,-1,-1). \]
By direct checking all possible case for $v_6'$ and $v_5'$, one concludes that 
it is impossible to get a reflexive $3$-polyhedron 
$P_9$ having exactly  $4$ faces of the type $F(4)$.

Thus,   vertices  $\{ v_1, \ldots, v_7 \}$ of $P$ can be always 
divided into $3$ pairs of centrally symmetric vectors.  
We assume that these centrally 
pairs are $\{ v_1, v_2 \},\; \{ v_3, v_4 \}, \;\{ v_5, v_6 \}$.  
Together with  $\{ v_8, v_9 \}$. Since there exist a $4$-dimensional 
cone in $\Sigma(P)$ which does not contain $v_7$, we can assume 
that $v_2, v_4, v_6, v_8$ is a basis of $\Z^4$. Now it is easy to 
check directly 
that $v_7 \in \{\pm v_2, \pm v_4, \pm v_6, \pm v_8\}$.

\hfill $\Box$

\begin{prop} 
Let $P$ be a Fano $4$-polyhedron. Then the following statements hold

{\rm (i)} $f_0(P) \leq 12$;

{\rm (ii)} if $f_0(P) = 12$, then $V(P)$ is the product of two 
Del Pezzo surfaces $S_3$ with the Picard rank $4$; 

{\rm (iii)}  if $f_0(P) = 11$, then $V(P)$ is the product of a   
Del Pezzo surface $S_3$ with the Picard rank $4$ and a Del Pezzo 
surface $S_2$ with the Picard rank $3$.
\label{11-12}
\end{prop}

\noindent
{\em Proof.} (i) By \ref{two.vert}, we can assume that there exists 
a primitive collection cointaining exactly $2$ vertices $\{ v_1, v_2\}$. 
Moreover, by \ref{symm-10}, we can assume that the corresponding primitive
relation is $v_1 + v_2 = v_0$. By \ref{except}, $D_0$ is a $\P^1$-bundle 
over a toric surface $S$. We know that  
the canonical class of $D_0$ is numerically 
effective (see \ref{num-eff}) and $P_0 = \pi_0(P)$ can contain only faces 
of the types $F(1)$ and $F(4)$ if $f_0(P) \geq 10$ (see 3.2-3.4). So the toric 
surface $S$ can have at most Picard number $4$, i.e., the valence of $D_0$ 
is at most $8$. Let $\{v_3, \ldots, v_k\}$ $(k \leq 8)$ be the set of 
vertices which are joined with $v_0$ by an edge of $P$. If one  
vertex $\pi_0(v_i) \in \{
\pi_0(v_3), \ldots, \pi_0(v_k) \}$ were a double vertex, 
then, by \ref{except}, 
$\Sigma_i(P)$ would have the combinatorial type $c_2$ (since the cone 
$\R_{\geq 0} \pi_0(v_i)$ is contained in exactly $4$ cones of $\Sigma_i(P)$). 
By \ref{c2-class} and \ref{c2-class2}, this would imply that $f_0(P) 
\leq 10$. If $\pi_0(v_1)$ and $\pi_0(v_2)$ are double points of $P_0$, 
then there exists at most $3$ vertices which are not not joined with 
$v_0$: at most $2$ vertices projecting to $\pi_0(v_1)$ and $\pi_0(v_2)$ and 
at most one centrally symmetric to $v_0$ vertex. This implies 
that $f_0(P) \leq 12$.    

(ii) It follows from the arguments in (i) that if $f_0(P) =12$, then 
$D_0$ is a $\P^1-bundle$ over a Del Pezzo surface $S_3$ obtained by 
blow ups of 
$3$-points in $\P^2$. Moreover, there exist primitive relations 
$v_0 + v_9 = v_1$, $v_0 + v_{10} = v_2$ and $v_0 + v_{11} =0$, i.e., $6$ 
vertices $v_0, v_1, v_2, v_9, v_{10}, v_{11}$ are contained in a 
$2$-dimensional linear subspace in $\R^4$ and their convex hull is isomorphic 
to a single Fano $2$-polyhedron with $6$ vertices. By \ref{un-comb}, the 
combinatorial type of $P$ is determined by $P_0$, in particular, 
$\{ v_3, v_4, v_5, v_6, v_7, v_8\}$ contains $9$ primitive collections 
consisting of $2$ vertices. Since, $P_0$ contains only faces of the type 
$F(1)$ and $F(4)$ the latter implies that
the convex hull of  $\{ v_3, v_4, v_5, v_6, v_7, v_8\}$ again must be 
isomorphic to a single Fano $2$-polyhedron with $6$ vertices. 

(iii) Let $f_0(P)=11$. Our previous arguments in (i) and (ii) show that 
we have a primitive relation $v_1 + v_2 = v_0$. It remains 
to consider the following two cases: 

{\sc Case 1}: Valence of $v_0$ equals $7$. Then there exist $3$ vertices  
which are not joined with 
$v_0$: $2$ vertices $v_8$, $v_9$ projecting to $\pi_0(v_1)$, 
$\pi_0(v_2)$ and one centrally symmetric to $v_0$ vertex $v_{10}$. As above 
we obtain that $v_0, v_1, v_2, v_8, v_{9}, v_{10}$ are contained in a 
$2$-dimensional linear subspace in $\R^4$ and their convex hull is isomorphic 
to a single Fano $2$-polyhedron with $6$ vertices and 
the convex hull of  $\{ v_3, v_4, v_5, v_6, v_7\}$ is 
isomorphic to a single Fano $2$-polyhedron with $5$ vertices. Therefore 
$V$ is isomorphic to $S_2 \times S_3$.

{\sc Case 2}: Valence of $v_0$ equals $8$. Then the vertices 
$v_0, v_1, v_2, v_9, v_{10}$ are contained in a 
$2$-dimensional linear subspace in $\R^4$ and their convex hull is isomorphic 
to a single Fano $2$-polyhedron with $5$ vertices. Without loss of generality 
we can assume that $v_0$ and $v_{10}$ are centrally symmetric (otherwise
we could replace $v_0$ by $v_1$ or by $v_2$). It remains to show 
that  the convex hull of $\{ v_3, v_4, v_5, v_6, v_7, v_8\}$ is isomorphic 
a single Fano $2$-polyhedron with $6$ vertices. The last property follows from 
that fact that the combinatorial type of $P$ is determined by 
$P_0$ (see \ref{un-comb}) and the faces of $P_0$ are of the type $F(1)$ or 
$F(4)$. 

\hfill $\Box$ 

\begin{prop} Let $P$ be a Fano $4$-polyhedron such that 
$f_0(P) \in \{9,10\}$ and 
all vertices of $P$ have valence $\geq 7$. Then the following 
statements hold

{\rm (i)} if $f_0(P) = 9$, then $P$ is isomorphic to the polyhedron 
from \ref{symm-10}(iii); 

{\rm (ii)} if $f_0(P) = 10$, then $P$ is isomorphic to the polyhedron 
from \ref{symm-10}(ii) or 
$V(P)$ is the product of two Del Pezzo surfaces $S_2$ with 
the Picard rank $3$.
\label{9-10}
\end{prop}
 
\noindent 
{\em Proof.} 
We can assume that there always exists a primitive relation 
of the type $v_1 + v_2 = v_3$. By \ref{except}, 
using  the contraction of the 
exceptional divisor $D_{3}$,   we obtain a smooth projective 
toric variety $V'$. We remark that $V$ must be a
 toric Fano $4$-fold coresponding 
to a Fano $4$-polyhedron $P'$, because vertices $v_1,v_2 \in P'$ 
have valence at least 
$7$, and the surface $D_{[1,2]}$ has the Picard number at least 
$3$ (see \ref{c-fano}).

To prove (i), we check all possibilities for the polyhedron $P'$ 
using the classification of Fano $4$-polyhedra with $8$ vertices. 
As a result, this completes the classification of Fano $4$-polyhedra with 
$9$ vertices. 

To prove (ii), we check all possibilities for the polyhedron $P'$ 
using the just received 
classification of Fano $4$-polyhedra with $9$ vertices.   
\hfill $\Box$
\bigskip

\begin{theo}  
There exist exactly 
$123$ different types of $4$-dimensional Fano polyhedra $P$. The maximal 
number of vertices of $P$ is $12$:

\medskip

\begin{center}
\begin{tabular}{|c|c|c|c|c|c|c|c|c|} \hline
\mbox{\em the number of vertices \/} & 5 & 6 & 7 & 
8 & 9 & \makebox{10} &  \makebox{11} & \makebox{12}  \\ \hline 
\mbox{\em the number of polyhedra\/} & 1 & 9 & 
\makebox{28} & \makebox{47} & \makebox{26} &  \makebox{10} & 
\makebox{1} & \makebox{1}  \\ \hline 
\end{tabular}
\end{center}
\end{theo}
\medskip

\newpage

\section{Table of toric Fano 4-folds} 

In the following table below we give the list  of  all  toric 
Fano 4-folds $V$ with their numerical characteristics: the products of 
characteristic classes $c_1^4$, $c_1^2 c_2$; the Betti numbers $b_2$, $b_4$; 
the dimensions $h^0 = h^0(V, -K_V)$ and $a(V):= dim\, Aut(V)$.  

\bigskip

\begin{tabular}{|c|c|c|c|c|c|c|l|} \hline 
    &      &     &     &     &    &     &   \\
 $n^0$  &  $c_1^4$  & $c_1^2 c_2$ &  $b_2$  &  $b_4$  & $a(V)$ & $h^0$  &
the type of Fano polytope 
\\    &      &     &     &     &    &     &   
\\ \hline 
  1 &  625 & 250 &  1  &  1  & 24 & 126 & ${\P}^4$
\\    &      &     &     &     &    &     &   
\\ \hline
  2 &  800 & 296 &  2  &  2  & 36 & 159 & $B_1, \; {\P}_{{\P}^3} 
({\cal O} \oplus {\cal O}(3))$
\\    &      &     &     &     &    &     &   
\\ \hline
  3 &  640 & 256 &  2  &  2  & 26 & 129 & $B_2, \; {\P}_{{\P}^3} 
({\cal O} \oplus {\cal O}(2))$
\\    &      &     &     &     &    &     &   
\\ \hline
  4 &  544 & 232 &  2  &  2  & 20 & 111 & $B_3, \; {\P}_{{\P}^3}
{\cal O} \oplus {\cal O}(1))$
\\    &      &     &     &     &    &     &   
\\ \hline
  5 &  512 & 224 &  2  &  2  & 18 & 105 & $B_4 , \; {\P}^1 \times {\P}^3$
\\    &      &     &     &     &    &     &   
\\ \hline
  6 &  512 & 224 &  2  &  2  & 18 & 105 & $B_5 , \; {\P}_{{\P}^1}
  ({\cal O} \oplus {\cal O} \oplus {\cal O} \oplus {\cal O}(1))$
\\    &      &     &     &     &    &     &   
\\ \hline
  7 &  594 & 240 &  2  &  3  & 24 & 120 & $C_1 , \; {\P}_{{\P}^2}
  ({\cal O} \oplus {\cal O} \oplus {\cal O}(2))$
\\    &      &     &     &     &    &     &   
\\ \hline
  8 &  513 & 222 &  2  &  3  & 18 & 105 & $C_2 , \; {\P}_{{\P}^2} 
  ({\cal O} \oplus {\cal O} \oplus {\cal O}(1))$
\\    &      &     &     &     &    &     &   
\\ \hline
  9 &  513 & 222 &  2  &  3  & 18 & 105 & $C_3 , \;{\P}_{{\P}^2} 
  ({\cal O} \oplus {\cal O}(1) \oplus {\cal O}(1))$
\\    &      &     &     &     &    &     &   
\\ \hline
 10 &  486 & 216 &  2  &  3  & 16 & 100 & $C_4 , \;{\P}^2 \times {\P}^2$
\\    &      &     &     &     &    &     &   
\\ \hline
 11 &  605 & 254 &  3  &  3  & 23 & 123 & $E_1$,  blow up of $\npt$  on $n^03$
\\    &      &     &     &     &    &     &   
\\ \hline
 12 &  489 & 222 &  3  &  3  & 17 & 101 & $E_2$, blow up of $\npt$ on $n^04$
\\    &      &     &     &     &    &     &   
\\ \hline
 13 &  431 & 206 &  3  &  3  & 14 & 90  & $E_3$, blow up of $\npt$ on $n^05$
\\    &      &     &     &     &    &     &   
\\ \hline
 14 &  592 & 244 &  3  &  4  & 24 & 120 & $D_1$, ${\npr}_{\npo \times \npt} 
 ({\cal O} \oplus {\cal O}(1,2))$
\\    &      &     &     &     &    &     &   
\\ \hline
 15 &  576 & 240 &  3  &  4  & 23 & 117 & $D_2$, $\npo$-bundle over 
$V({\cal B}_1)$
\\    &      &     &     &     &    &     &   
\\ \hline
 16 &  560 & 236 &  3  &  4  & 22 & 114 & $D_3$, $\npo$-bundle over
$V({\cal B}_2)$
\\    &      &     &     &     &    &     &   
\\ \hline
 17 &  560 & 236 &  3  &  4  & 22 & 114 & $D_4$, $\npo$-bundle over 
$V({\cal B}_3)$
\\    &      &     &     &     &    &     &   
\\ \hline 
\end{tabular}

\newpage

\begin{tabular}{|c|c|c|c|c|c|c|l|} \hline 
    &      &     &     &     &    &     &   \\
 $n^0$  &  $c_1^4$  & $c_1^2 c_2$ &  $b_2$  &  $b_4$  & $a(V)$ & $h^0$  &
the type of Fano polytope 
\\    &      &     &     &     &    &     &   
\\ \hline
 18 &  496 & 220 &  3  &  4  & 18 & 102 & $D_5$, $\npo \times {\npr}_{\npt}
 ({\cal O} \oplus {\cal O}(2))$
\\    &      &     &     &     &    &     &   
\\ \hline
 19 &  496 & 220 &  3  &  4  & 18 & 102 & $D_6$, ${\npr}_{\npo \times \npt} 
 ({\cal O} \oplus {\cal O}(1,1))$
\\    &      &     &     &     &    &     &   
\\ \hline
 20 &  486 & 216 &  3  &  4  & 18 & 100 & $D_7$, ${\npr}_{\npo \times \npo} 
  ({\cal O} \oplus {\cal O} \oplus {\cal O}(1,1))$
\\    &      &     &     &     &    &     &   
\\ \hline
21 &  432 & 204 &  3  &  4  & 14 &  90 & $D_8$, $\npo$-bundle over 
$V({\cal B}_2)$
\\    &      &     &     &     &    &     &   
\\ \hline
 22 &  464 & 212 &  3  &  4  & 16 &  96 & $D_9$, $\npo$-bundle over 
$V({\cal B}_2)$
\\    &      &     &     &     &    &     &   
\\ \hline
 23 &  464 & 212 &  3  &  4  & 16 &  96 & $D_{10}$, $\npo$-bundle over 
$V({\cal B}_3)$
\\    &      &     &     &     &    &     &   
\\ \hline
 24 &  459 & 210 &  3  &  4  & 15 &  95 & $D_{11}$, ${\npr}_{F_1} ({\cal O} 
 \oplus {\cal O} \oplus {\cal O}(l))$ \\
    &      &     &     &     &    &     & $l$ is a curve of index 1 on $F_1$
\\ \hline
 25 &  448 & 208 &  3  &  4  & 15 & 93 & $D_{12}$, $\npo \times {\npr}_
{\npt} ({\cal O} \oplus {\cal O}(1))$
\\    &      &     &     &     &    &     &   
\\ \hline
 26 &  432 & 204 &  3  &  4  & 14 & 90 & $D_{13}$, $\npo \times \npo 
\times \npt$
\\    &      &     &     &     &    &     &   
\\ \hline
 27 &  432 & 204 &  3  &  4  & 14 & 90 & $D_{14}$, $\npo \times 
\times {\npr}_{\npo} ({\cal O} \oplus {\cal O} \oplus {\cal O}(1))$
\\    &      &     &     &     &    &     &   
\\ \hline
 28 &  432 & 204 &  3  &  4  & 14 & 90 & $D_{15}$, $F^1 \times \npt$
\\    &      &     &     &     &    &     &   
\\ \hline
 29 &  432 & 204 &  3  &  4  & 14 & 90 & $D_{16}$, $\npo$-bundle over 
$V({\cal B}_1)$
\\    &      &     &     &     &    &     &   
\\ \hline
 30 &  405 & 198 &  3  &  4  & 12 & 85 & $D_{17}$, ${\npr}_{\npo \times \npo} 
({\cal O} \oplus {\cal O}(1,0) \oplus {\cal O}(0,1))$
\\    &      &     &     &     &    &     &   
\\ \hline
 31 &  400 & 196 &  3  &  4  & 12 & 84 & $D_{18}$, ${\npr}_{\npo \times \npt} 
({\cal O} \oplus {\cal O}(-1,2))$
\\    &      &     &     &     &    &     &   
\\ \hline
 32 &  400 & 196 &  3  &  4  & 12 & 84 & $D_{19}$, ${\npr}_{\npo \times \npt} 
({\cal O} \oplus {\cal O}(-1,1))$
\\    &      &     &     &     &    &     &   
\\ \hline
 33 &  529 & 226 &  3  &  5  & 20 & 108 & $G_1$, is not contractible smoothly
\\   &      &     &     &     &    &     &   
\\ \hline
 34 &  450 & 204 &  3  &  5  & 16 & 93 & $G_2$, 
blow up of a surface 
\\  &   &  &    &    &  &  & 
$\; \; \; \; \;$ 
$\npop$ on $n^0 8$
\\ \hline
 35 &  433 & 202 &  3  &  5  & 14 & 90 & $G_3$, 
blow up of a curve 
\\  &   &  &    &    &  &  & 
$\; \; \; \; \;$ 
$\npo$ on $n^0 9$
\\ \hline
 36 &  417 & 198 &  3  &  5  & 13 & 87 & $G_4$, 
blow up of a surface 
\\  &   &  &    &    &  &  & 
$\; \; \; \; \;$ 
$F_1 $ on $n^0 8$
\\ \hline
 37 &  406 & 196 &  3  &  5  & 12 & 85 & $G_5$,  
blow up of a surface 
\\  &   &  &    &    &  &  & 
$\; \; \; \; \;$ 
$\npop$ on $n^0 9$
\\ \hline
\end{tabular}

\newpage

\begin{tabular}{|c|c|c|c|c|c|c|l|} \hline 
    &      &     &     &     &    &     &   \\
 $n^0$  &  $c_1^4$  & $c_1^2 c_2$ &  $b_2$  &  $b_4$  & $a(V)$ & $h^0$  &
the type of Fano polytope 
\\    &      &     &     &     &    &     &   
\\ \hline
 38 &  401 & 194 &  3  &  5  & 12 & 84 & $G_6$,  
blow up of a surface 
\\  &   &  &    &    &  &  & 
$\; \; \; \; \;$ 
$\npop$ on $n^0 10$
\\ \hline
 39 &  588 & 240 &  4  &  5  & 22 & 114 & $H_1$, 
blow up of two surfaces 
\\  &   &  &    &    &  &  & 
$\; \; \; \; \;$ 
$\npt$ on $n^0 3$
\\ \hline
 40 &  505 & 226 &  4  &  5  & 19 & 104 & $H_2$,  
blow up of a surface 
\\  &   &  &    &    &  &  & 
$\; \; \; \; \;$ 
$\npt$ on $n^0 16$
\\ \hline 
 41 &  478 & 220 &  4  &  5  & 17 &  99 & $H_3$,  
blow up of a surface 
\\  &   &  &    &    &  &  & 
$\; \; \; \; \;$ 
$\npt$ on $n^0 14$
\\ \hline
 42 &  447 & 210 &  4  &  5  & 16 &  93 & $H_4$,  
blow up of two surfaces 
\\  &   &  &    &    &  &  & 
$\; \; \; \; \;$ 
$\npt$ on $n^0 8$
\\ \hline
 43 &  415 & 202 &  4  &  5  & 14 &  87 & $H_5$,  
blow up of a surface 
\\  &   &  &    &    &  &  & 
$\; \; \; \; \;$ 
$\npt$ on $n^0 19$
\\ \hline
 44 &  409 & 202 &  4  &  5  & 13 &  96 & $H_6$,  
blow up of a surface 
\\  &   &  &    &    &  &  & 
$\; \; \; \; \;$ 
$\npt$ on $n^0 22$
\\ \hline
 45 &  382 & 196 &  4  &  5  & 11 &  81 & $H_7$,  
blow up of two surfaces 
\\  &   &  &    &    &  &  & 
$\; \; \; \; \;$ 
$\npt$ on $n^0 3$
\\ \hline
 46 &  378 & 192 &  4  &  5  & 12 &  80 & $H_8$,  $\npt \times S_2$
\\  &   &  &    &    &  &  & 
\\ \hline
 47 &  367 & 190 &  4  &  5  & 11 &  78 & $H_9$, 
blow up of two surfaces 
\\  &   &  &    &    &  &  & 
$\; \; \; \; \;$ 
$\npt$ on $n^0 8$
\\ \hline
 48 &  351 & 186 &  4  &  5  & 10 &  75 & $H_{10}$, 
blow up of a surface 
\\  &   &  &    &    &  &  & 
$\; \; \; \; \;$ 
$\npt$ on $n^0 19$
\\ \hline
 49 &  480 & 216 &  4  &  6  & 18 &  99 & $L_1$, ${\npr}_{\npo \times 
\npo \times \npo } ({\cal O} \oplus {\cal O}(1,1,1))$
\\  &   &  &    &    &  &  & 
\\ \hline
 50 &  464 & 212 &  4  &  6  & 17 &  96 & $L_2$, $\npo$-bundle over 
$V({\cal C}_1)$
\\  &   &  &    &    &  &  & 
\\ \hline
 51 &  448 & 208 & 4   &  6  & 16 &  93  &$L_3$, 
${\npr}_{\npo \times F_1 } ({\cal O} \oplus {\cal O}(1) \otimes {\cal O}(l)$
\\  &   &  &    &    &  &  & 
$\; \; \; \; \;$ 
$l$ is a curve of index 1 on $F_1$
\\ \hline
 52 &  432 & 204 &  4  &  6  & 15 &  90 & $L_4$, $\npo$-bundle over 
$V({\cal C}_3)$
\\  &   &  &    &    &  &  & 
\\ \hline
 53 &  416 & 200 &  4  &  6  & 14 &  87 & $L_5$, $\npo \times 
{\npr}_{\npo \times \npo} ({\cal O} \oplus {\cal O}(1,1))$
\\  &   &  &    &    &  &  & 
\\ \hline
 54 &  400 & 196 &  4  &  6  & 13 &  84 &$L_6$, 
${\npr}_{\npo \times F_1 } ({\cal O} \oplus  {\cal O}(l))$
\\  &   &  &    &    &  &  & 
$\; \; \; \; \;$ 
$l$ is a curve of index 1 on $F_1$
\\ \hline
 55 &  384 & 192 &  4  &  6  & 12 &  81 &$L_7$, $F_1 \times F_1$
\\  &   &  &    &    &  &  & 
\\ \hline
 56 &  384 & 192 &  4  &  6  & 12 &  81 &$L_8$, $\npo \times \npo \times 
\npo \times \npo$
\\  &   &  &    &    &  &  & 
\\ \hline
 57 &  384 & 192 &  4  &  6  & 12 &  81 &$L_9$, $\npo \times \npo \times 
F_1$
\\  &   &  &    &    &  &  & 
\\ \hline
\end{tabular}

\newpage

\begin{tabular}{|c|c|c|c|c|c|c|l|} \hline
    &      &     &     &     &    &     &   \\
 $n^0$  &  $c_1^4$  & $c_1^2 c_2$ &  $b_2$  &  $b_4$  & $a(V)$ & $h^0$  &
the type of Fano polytope 
\\    &      &     &     &     &    &     &   
\\ \hline
 58 &  384 & 192 &  4  &  6  & 12 &  81 &$L_{10}$, $\npo$-bundle over 
$V({\cal D}_2)$
\\  &   &  &    &    &  &  & 
\\ \hline
 59 &  352 & 184 &  4  &  6  & 10 &  75 &$L_{11}$, $\npo \times 
V({\cal D}_2)$
\\  &   &  &    &    &  &  & 
\\ \hline
 60 &  352 & 184 &  4  &  6  & 10 &  75 &$L_{12}$, 
${\npr}_{\npo \times F_1 } ({\cal O} \oplus {\cal O}(1) \otimes {\cal O}(-l))$
\\  &   &  &    &    &  &  & 
$\; \; \; \; \;$ 
$l$ is a curve of index 1 on $F_1$
\\ \hline 
61 & 352  & 184 &  4  &  6  & 10 & 75 & $L_{13}$, ${\npr}_{\npo \times 
\npo \times \npo} ({\cal O} \oplus {\cal O} (1,1,-1))$
\\  &      &     &     &     &    &    & 
$\; \; \; \; \; $ 
\\ \hline
 62 & 496  & 220 &  4  &  6  & 19 & 102 & $I_1$, 
blow up of a surface 
\\  &      &     &     &     &    &    & 
$\; \; \; \; \; $ 
$\npop$ on $n^0 17$
\\ \hline
 63 & 463  & 214 &  4  &  6  & 16 & 96 & $I_2$, 
blow up of a surface 
\\  &      &     &     &     &    &    & 
$\; \; \; \; \; $ 
$\npop$ on $n^0 19$
\\ \hline
 64 & 442  & 208 &  4  &  6  & 15 & 92 & $I_3$, 
blow up of a surface 
\\  &      &     &     &     &    &    & 
$\; \; \; \; \; $ 
$\npop$ on $n^0 21$
\\ \hline
 65 & 433  & 202 &  4  &  6  & 13 & 87 & $I_4$, 
blow up of a surface 
\\  &      &     &     &     &    &    & 
$\; \; \; \; \; $ 
$\npop$ on $n^0 23$
\\ \hline
 66 & 415  & 202 &  4  &  6  & 13 & 87 & $I_5$, 
blow up of a surface 
\\  &      &     &     &     &    &    & 
$\; \; \; \; \; $ 
$\npt$ on $n^0 23$
\\ \hline
 67 & 411  & 198 &  4  &  6  & 14 & 86 & $I_6$, 
blow up of a surface 
\\  &      &     &     &     &    &    & 
$\; \; \; \; \; $ 
$\npop$ on $n^0 23$
\\ \hline
 68 & 400  & 196 &  4  &  6  & 13 & 84 & $I_7$, $\npo \times V({\cal D}_1)$
\\  &      &     &     &     &    &    & 
\\ \hline
 69 & 384  & 192 &  4  &  6  & 12 & 81 & $I_8$, 
blow up of a surface 
\\  &      &     &     &     &    &    & 
$\; \; \; \; \; $ 
$\npop$ on $n^0 21$
\\ \hline
 70 &  390 & 192 &  4  &  6  & 13 & 82 & $I_9$,  
blow up of a surface 
\\  &      &     &     &     &    &    & 
$\; \; \; \; \; $ 
$\npop$ on $n^0 25$
\\ \hline
 71 &  389 & 194 &  4  &  6  & 12 & 82 & $I_{10}$,   
blow up of a surface 
\\  &      &     &     &     &    &    & 
$\; \; \; \; \; $ 
$\npop$ on $n^0 28$
\\ \hline
 72 &  384 & 192 &  4  &  6  & 12 & 81 & $I_{11}$, $\npo$-bundle over 
$V({\cal D}_1)$   \\  &      &     &     &     &    &    & 
\\ \hline
 73 &  347 & 182 &  4  &  6  & 10 & 74 & $I_{12}$,   
blow up of a surface 
\\  &      &     &     &     &    &    & 
$\; \; \; \; \; $ 
$\npop$ on $n^0 28$
\\ \hline
 74 &  368 & 188 &  4  &  6  & 11 & 78 & $I_{13}$, $\npo \times 
V({\cal D}_2)$        \\  &      &     &     &     &    &    & 
\\ \hline
 75 &  357 & 186 &  4  &  6  & 10 & 76 & $I_{14}$,  
blow up of a surface 
\\  &      &     &     &     &    &    & 
$\; \; \; \; \; $ 
$F_1$ on $n^0 27$
\\ \hline
 76 &  337 & 178 &  4  &  6  & 10 & 72 & $I_{15}$,     
blow up of a surface 
\\  &      &     &     &     &    &    & 
$\; \; \; \; \; $ 
$\npop$ on $n^0 32$
\\ \hline
 77 &  385 & 190 &  4  &  7  & 12 & 81 & $M_1$, is not contractible smoothly
\\  &      &     &     &     &    &    & 
$\; \; \; \; \; $ 
\\ \hline
\end{tabular}

\newpage

\begin{tabular}{|c|c|c|c|c|c|c|l|} \hline
    &      &     &     &     &    &     &   \\
 $n^0$  &  $c_1^4$  & $c_1^2 c_2$ &  $b_2$  &  $b_4$  & $a(V)$ & $h^0$  &
the type of Fano polytope 
\\    &      &     &     &     &    &     &   
\\ \hline
 78 &  417 & 198 &  4  &  7  & 14 & 87 & $M_2$, is not contractible smoothly
\\  &      &     &     &     &    &    & 
$\; \; \; \; \; $ 
\\ \hline
 79 &  369 & 186 &  4  &  7  & 11 & 78 & $M_3$,   
blow up of a surface 
\\  &      &     &     &     &    &    & 
$\; \; \; \; \; $ 
$F_1$ on $n^0 35$
\\ \hline
 80 &  369 & 186 &  4  &  7  & 11 & 78 & $M_4$,    
blow up of a surface 
\\  &      &     &     &     &    &    & 
$\; \; \; \; \; $ 
$\npop$ on $n^0 35$
\\ \hline 
81 &  364 & 184 &  4  &  7  & 11 & 77 &  $M_5$,  
blow up of a surface 
$\npop$ 
\\  &      &     &     &     &    &    & 
$\; \; \; \; \; $ 
on $\npt \times \npt$,
and then blow up of $F_1$
\\ \hline
 82 &  368 & 188 &  4  &  7  & 10 & 78 & $J_1$,   
blow up of a surface 
\\  &      &     &     &     &    &    & 
$\; \; \; \; \; $ 
$\npop$ on $n^0 35$
\\ \hline
 83 &  326 & 176 &  4  &  7  &  8 & 70 & $J_2$,      
blow up of a surface 
\\  &      &     &     &     &    &    & 
$\; \; \; \; \; $ 
$\npop$ on $n^0 35$
\\ \hline
 84 &  422 & 208 &  5  &  8  & 16 & 92 & $Q_1$,   
blow up of two surfaces 
\\  &      &     &     &     &    &    & 
$\; \; \; \; \; $ 
$\npop$ on $n^0 20$
\\ \hline
 85 &  405 & 198 &  5  &  8  & 14 & 85 & $Q_2$,    
blow up of two surfaces 
\\  &      &     &     &     &    &    & 
$\; \; \; \; \; $ 
$F_1$ on $n^0 28$
\\ \hline
 86 &  394 & 196 &  5  &  8  & 13 & 83 & $Q_3$,    
blow up of a surface 
\\  &      &     &     &     &    &    & 
$\; \; \; \; \; $ 
$\npop$ on $n^0 49$
\\ \hline
 87 &  405 & 198 &  5  &  8  & 14 & 85 & $Q_4$,   
blow up of a surface 
\\  &      &     &     &     &    &    & 
$\; \; \; \; \; $ 
$\npop$ on $n^0 51$
\\ \hline
 88 &  373 & 190 &  5  &  8  & 12 & 79 & $Q_5$,    
blow up of a surface 
\\  &      &     &     &     &    &    & 
$\; \; \; \; \; $ 
$F_1$ on $n^0 51$
\\ \hline
 89 &  368 & 188 &  5  &  8  & 12 & 78 & $Q_6$,
     $\npo \times V({\cal E}_1)$            
     \\  &      &     &     &     &    &    & 
$\; \; \; \; \; $ 
\\ \hline
 90 &  363 & 186 &  5  &  8  & 12 & 77 & $Q_7$,     
blow up of a surface 
\\  &      &     &     &     &    &    & 
$\; \; \; \; \; $ 
$\npop$ on $n^0 55$
\\ \hline
 91 &  352 & 184 &  5  &  8  & 11 & 75 & $Q_8$,   
 $\npo \times  V({\cal E}_2)$         
\\  &      &     &     &     &    &    & 
$\; \; \; \; \; $ 
\\ \hline
 92 &  341 & 182 &  5  &  8  & 10 & 73 & $Q_9$,      
blow up of a surface 
\\  &      &     &     &     &    &    & 
$\; \; \; \; \; $ 
$\npop$ on $n^0 51$
\\ \hline
 93 &  336 & 180 &  5  &  8  & 10 & 72 & $Q_{10}$,  $F_1 \times S_2$   
\\  &      &     &     &     &    &    & 
$\; \; \; \; \; $ 
\\ \hline
 94 &  336 & 180 &  5  &  8  & 10 & 72 & $Q_{11}$,  $\npop \times S_2$  
\\  &      &     &     &     &    &    & 
\\ \hline
 95 &  331 & 188 &  5  &  8  & 10 & 71 & $Q_{12}$,  
blow up of two surfaces 
\\  &      &     &     &     &    &    & 
$\; \; \; \; \; $ 
$\npop$ on $n^0 30$
\\ \hline
 96 &  330 & 180 &  5  &  8  &  9 & 71 & $Q_{13}$,    
blow up of two surfaces 
\\  &      &     &     &     &    &    & 
$\; \; \; \; \; $ 
$\npop$ on $n^0 20$
\\ \hline
 97 &  325 & 178 &  5  &  8  &  9 & 70 & $Q_{14}$,      
blow up of two surfaces
\\  &      &     &     &     &    &    & 
$\; \; \; \; \; $ 
$F_1$ on $n^0 24$
\\ \hline
\end{tabular}

\newpage

\begin{tabular}{|c|c|c|c|c|c|c|l|} \hline
    &      &     &     &     &    &     &   \\
 $n^0$  &  $c_1^4$  & $c_1^2 c_2$ &  $b_2$  &  $b_4$  & $a(V)$ & $h^0$  &
the type of Fano polytope 
\\    &      &     &     &     &    &     &   
\\ \hline
 98 &  320 & 176 &  5  &  8  &  9 & 69 & $Q_{15}$,  
 $\npo \times V({\cal E}_4)$        
\\  &      &     &     &     &    &    & 
\\ \hline
 99 &  310 & 172 &  5  &  8  &  9 & 67 & $Q_{16}$,        
blow up of a surface 
\\  &      &     &     &     &    &    & 
$\; \; \; \; \; $ 
$\npop$ on $n^0 61$
\\ \hline
 100 &  299 & 170 &  5  &  8  &  8 & 65 & $Q_{17}$,        
blow up of a surface 
\\  &      &     &     &     &    &    & 
$\; \; \; \; \; $ 
$\npop$ on $n^0 55$
\\ \hline 
101 &  364 & 196 &  5  &  6  & 10 & 78 & $K_1$,        
blow up of three surfaces 
\\  &      &     &     &     &    &    & 
$\; \; \; \; \; $ 
$\npt$ on $n^0 7$
\\ \hline
 102 &  354 & 192 &  5  &  6  & 10 & 76 & $K_2$,   
blow up of two surfaces 
\\  &      &     &     &     &    &    & 
$\; \; \; \; \; $ 
$\npt$ on $n^0 22$
\\ \hline
 103 &  334 & 184 &  5  &  6  & 10 & 72 & $K_3$,    
blow up of three surfaces 
\\  &      &     &     &     &    &    & 
$\; \; \; \; \; $ 
$\npt$ on $n^0 8$
\\ \hline
 104 & 324 & 180 &  5  &  6  & 10 & 70 & $K_4$,   $\npt \times S_3$    
\\  &      &     &     &     &    &    & 
\\ \hline
 105 &  332 & 176 &  5  &  9  & 10 & 71 & $R_1$,     
blow up of a surface 
\\  &      &     &     &     &    &    & 
$\; \; \; \; \; $ 
$\npop$ on $n^0 78$
\\ \hline
 106 &  321 & 174 &  5  &  9  &  9 & 69 & $R_2$,       
blow up of a surface 
\\  &      &     &     &     &    &    & 
$\; \; \; \; \; $ 
$S_1$ on $n^0 78$
\\ \hline
 107 &  305 & 170 &  5  &  9  &  8 & 66 & $R_3$,        
blow up of a surface 
\\  &      &     &     &     &    &    & 
$\; \; \; \; \; $ 
$\npop$ on $n^0 79$
\\ \hline
 108 &  310 & 172 &  5  &  9  &  8 & 67 &     
$\; \; \; \; \; $ 
blow up of a surface 
\\  &      &     &     &     &    &    & 
$\; \; \; \; \; $ 
$S_2$ on $\npo \times V({\cal D}_{2})$ see 3.4.1
\\ \hline
 109 &  308 & 176 &  6  & 10  &  8 & 67 & $U_1$,          
blow up of three surfaces 
\\  &      &     &     &     &    &    & 
$\; \; \; \; \; $ 
$\npop$ on $n^0 20$
\\ \hline
 110 &  298 & 172 &  6  & 10  &  8 & 65 & $U_2$,          
blow up of three surfaces 
\\  &      &     &     &     &    &    & 
$\; \; \; \; \; $ 
$F_1$ on $n^0 24$
\\ \hline
 111 &  298 & 172 &  6  & 10  &  8 & 65 & $U_3$,        
blow up of two surfaces 
\\  &      &     &     &     &    &    & 
$\; \; \; \; \; $ 
$\npop$ on $S_1 \times S_1$
\\ \hline
 112 &  288 & 168 &  6  & 10  &  8 & 63 & $U_4$,    $S_1 \times  S_3$  
\\  &      &     &     &     &    &    & 
\\ \hline
 113 &  288 & 168 &  6  & 10  &  8 & 63 & $U_5$,  $\npop \times S_3$  
\\  &      &     &     &     &    &    & 
\\ \hline
 114 &  288 & 168 &  6  & 10  &  8 & 63 & $U_6$,   $\npo \times 
V({\cal F}_{1})$ 
\\  &      &     &     &     &    &    & 
\\ \hline
 115 &  278 & 164 &  6  & 10  &  8 & 61 & $U_7$,    
blow up of three surfaces 
\\  &      &     &     &     &    &    & 
$\; \; \; \; \; $ 
$\npop$ on $n^0 30$
\\ \hline
 116 &  268 & 160 &  6  & 10  &  8 & 59 & $U_8$,   
blow up of two surfaces 
\\  &      &     &     &     &    &    & 
$\; \; \; \; \; $ 
$\npop$ on $\P^1 \times V({\cal C}_5)$  
\\ \hline
 117 &  307 & 166 &  5  & 11  &  8 & 66 &   see \ref{symm-10}(ii)
\\  &      &     &     &     &    &    & 
\\ \hline
\end{tabular}
\newpage

\begin{tabular}{|c|c|c|c|c|c|c|l|} \hline
    &      &     &     &     &    &     &   \\
 $n^0$  &  $c_1^4$  & $c_1^2 c_2$ &  $b_2$  &  $b_4$  & $a(V)$ & $h^0$  &
the type of Fano polytope 
\\    &      &     &     &     &    &     &   
\\ \hline
 118 &  230 & 140 &  6  & 16  &  4 & 51 &   see \ref{symm-10}(iii)  
\\  &      &     &     &     &    &    & 
\\ \hline
 119 &  294 & 168 &  6  & 11  &  8 & 64 &  $S_2 \times S_2$  
\\  &      &     &     &     &    &    & 
\\ \hline
 120 &  252 & 156 &  7  & 14  &  6 & 56 &  $S_2 \times S_3$     
\\  &      &     &     &     &    &    & 
\\ \hline 

121 &  216 & 144 &  8  & 18  &  4 & 49 &  $S_3 \times S_3$  
\\  &      &     &     &     &    &    & 

\\ \hline 

122 &  322  & 172 &  4  &  8  &  8 & 69 &  $Z_1$, blow up of a surface                          
\\  &      &     &     &     &    &    & 
$\; \; \; \; \; $ 
$S_2$ on $n^0 38$
\\ \hline 

123 & 327  & 174 &  4  &  8  & 8  & 70 &  $Z_2$, blow up of a surface                          
\\  &      &     &     &     &    &    & 
$\; \; \; \; \; $ 
$S_2$ on $n^0 36$
\\ \hline
\end{tabular}
\bigskip 
\bigskip

Using  the generalized Euler exact sequence for 
smooth projective toric varieties \cite{bat.mel}, one obtains the following 
formula for the dimension of the group of biregular authomorphisms 
of a toric Fano $4$-fold $V$: 
\[ a(V) = 22 + 3 b_2 - b_4 + \frac{2c_1^4 - 5c_1 c_2}{12}. \]
If $W \subset V$ is a smooth Calabi-Yau $3$-fold, then the Hodge 
number $h^{2,1}(W)$ is determined by the formula
\[ h^{2,1}(W) =  h^0(V, -K_V) -a(V) -1 = \frac{c_1^2c_2}{2} + b_4
- 3b_2 - 22, \]
where $h^0(V, -K_V) = 1 + \frac{1}{12}( c_1^2c_2 + 2c_1^4)$.


\begin{thebibliography}{99}                            

\bibitem{bat1} Batyrev, V.V.: {\em Toroidal Fano 3-folds}, 
Math. USSR-Izv. 19, 13-25 (1982)
    
\bibitem{bat2}  Batyrev, V.V.: {\em Boundness of the degree of  
higher-dimensional  toric Fano varieties}, Moscow Univ. Math. Bull. 
37, 28-33 (1982)


\bibitem{bat4}  Batyrev, V.V.: {\em Higher-dimensional toric varieties 
with ample anticanonical class}. PhD Thesis (in Russian), Moscow State 
University, 1984. 

\bibitem{bat.mel} Batyrev, V.V.; Mel'nikov, D.A.:
{\em A theorem on non-extensibility of toric varieties},  
Mosc. Univ. Math. Bull. {\bf 41}, No.3, 23-27 (1986)


\bibitem{bat5}  Batyrev, V.V.: 
{\em On the classification of smooth projective toric varieties}, 
Tohoku Math. J., II. Ser. 43, No.4, 569-585 (1991)


\bibitem{bat.dual} Batyrev, V. V.:
{\em Dual polyhedra and mirror symmetry for Calabi-Yau hypersurfaces in 
toric varieties}, J. Algebr. Geom. 3, No.3, 493-535 (1994)

\bibitem{bat.cox} Batyrev, V. V.; Cox, D. A.: 
{\em On the Hodge structure of projective hypersurfaces in toric 
varieties}, Duke Math. J. 75, No.2, 293-338 (1994)

\bibitem{borel} Borel, A.: 
{\em Linear algebraic groups},  2nd ed., Graduate Texts in Mathematics, 
{\bf 126}, New York, Springer-Verlag (1991)

\bibitem{borisovs} Borisov, A.A.; Borisov, L.A.:
{\em Singular toric Fano varieties}, Math. USSR, Sb. 75, No.1, 
277-283 (1993)

\bibitem{cox} Cox, D. A.: 
{\em The homogeneous coordinate ring of a toric variety}, 
J. Algebr. Geom. {\bf 4}, No.1, 17-50 (1995)

\bibitem{cutkosky} Cutkosky, S. D.: 
{\em On Fano 3-folds},  
Manuscr. Math. 64, No.2, 189-204 (1989). 


\bibitem{dan2} Danilov V.I.: {\em The geometry of toric varieties},  
Uspekhi  Math. Nauk (2), 33, 85-134(1978); Math. URSS - Izv., 17,97-154 (1981)

\bibitem{evertz} Evertz, S.: {\em Zur Klassifikation 4-dimensionaler 
Fano-Variet\"{a}ten}. Diplomarbeit (in German), Math. Inst. der 
Ruhr-Universit\"{a}t 
Bochum, September 1988. 

\bibitem{ewald} Ewald, G.: {\em On  the  classification  of  toric  
Fano  varieties}.  Discrete Comput. Geom. 3, 49-54 (1988)

\bibitem{fulton} Fulton, W.: {\em 
Introduction to toric varieties}, The 1989 William H. Roever lectures 
in geometry, Annals of Mathematics Studies {\bf 131},  Princeton University 
Press, (1993) 

\bibitem{grun} Gr\"unbaum, B.; Streedharan, V.P.: {\em An enumeration of 
simplicial  4-polytopes with 8  vertices},  J.  
Combinatorial  Theory,  2,  435-465 (1967)

\bibitem{isk1} Iskovskih, V.A.: {\em Fano 3-folds I},  
Math.  USSR-Izv,  485-527 (1977) 

\bibitem{isk2} Iskovskih, V.A.: {\em Fano 3-folds II}, Math. USSR - 
Izv., 12, 469-506 (1978)

\bibitem{kl} Kleinschmidt, P.: {\em A classification of toric varieties 
with few generators}, 
Aequationes Math. 35, No.2/3, 254-266 (1988)

\bibitem{KMM1} Kollar, J.; Miyaoka, Y.; Mori, S.: 
{\em Rational connectedness and boundedness of Fano manifolds}, 
J. Differ. Geom. 36, No.3, 765-779 (1992).

\bibitem{KMM2} Kollar, J.; Miyaoka, Y.; Mori, S.: 
{\em Rational curves on Fano varieties}, 
in {\sl Classification of irregular varieties, minimal models 
and abelian varieties}, Proc. Conf., Trento/Italy 1990,
Lect. Notes Math. {\bf 1515}, 100-105 (1992).

\bibitem{mabuchi} 
Mabuchi, T.: {\em Einstein - K\"ahler  forms,  Futaki  invariants  and 
convex geometry on  toric  Fano  varieties},  
Osaka  J.  Math,  24, 705-737 (1987)

\bibitem{manin} Manin, Yu. I.: {\em Cubic forms. Algebra, geometry, 
arithmetic}, 2nd ed., Amsterdam-New York-Oxford: North-Holland, (1986)

\bibitem{mori1} Mori, S.: 
{\em Threefolds whose canonical bundles are not numerically effective}, 
Ann. Math., II. Ser. 116, 133-176 (1982) 

\bibitem{mori2} Mori, S.: 
{\em Cone of curves, and Fano 3-folds}, 
Proc. Int. Congr. Math., Warszawa 1983, Vol. 1, 747-752 (1984)

\bibitem{MM1} Mori, S.;  Mukai,  S.: {\em Classification  of  
Fano  3-folds  with $B_2 >2$}, Manusripta Math. 36, 147-162 (1981)

\bibitem{MM2}  Mori, S.; Mukai,  S.: 
{\em On Fano 3-folds with  $B_2 > 2$}, in {\sl Algebraic varieties and 
analytic varieties}, Proc. Symp., Tokyo 1981, Adv. Stud. Pure Math. 1, 
101-129 (1983)  

\bibitem{MM3} 
Mori, S.; Mukai, S.: {\em Classification of Fano $3$-folds with 
$B\sb 2\geq 2$. I}, in  {\sl  Algebraic and topological theories} 
Nagata, M. (ed.) et al.,  
Papers from the symposium dedicated to the
memory of Dr. Takehiko Miyata held in Kinosaki, 
October 30-November 9, 1984. Tokyo: Kinokuniya Company
Ltd. 496-545 (1986)

\bibitem{murre} Murre, J.P.: 
{\em Classification of Fano threefolds according to Fano and 
Iskovskih}, in 
{\sl Algebraic threefolds}, Proc. 2nd 1981 Sess. C.I.M.E., 
Varenna/Italy 1981, Lect. Notes Math. 947, 35-92
(1982).


\bibitem{nadel} Nadel, A. M.: 
{\em The boundedness of degree of Fano varieties with Picard number one}, 
J. Am. Math. Soc. 4, No.4, 681-692 (1991). 

\bibitem{nakagawa0} Nakagawa, Y.: {\em A letter to V. Batyrev}, 
November 13, 1991. 

\bibitem{nakagawa1} Nakagawa, Y.: {\em 
Einstein-Kaehler toric Fano fourfolds}, 
Tohoku Math. J., II. Ser. 45, No.2, 297-310 (1993). 

\bibitem{nakagawa2} Nakagawa, Y.:
{\em Classification of Einstein-Kaehler toric Fano fourfolds}, 
Tohoku Math. J., II. Ser. 46, No.1, 125-133 (1994). 

\bibitem{oda} Oda, T.: {\em Convex bodies and algebraic geometry.
An introduction to the theory of toric varieties}, 
Ergebnisse der Mathematik und ihrer Grenzgebiete. 3. Folge, Bd. 15. 
Berlin, Springer-Verlag (1988). 

\bibitem{reid} Reid, M.: {\em Decomposition of toric morphisms},  
in {\sl  Arithmetic and geometry}, Vol. II: {\sl Geometry}, 
Prog. Math. 36, 395-418
(1983)

\bibitem{vosk.klyachko} Voskresenskij, V.E.; Klyachko A.A.: 
{\em Toroidal Fano varieties  and root systems}, Math. USSR-Izv. 
24, 221-244 (1985) 

\bibitem{wat} Watanabe, K.; Watanabe, M: {\em The  classification  
of  Fano  3-folds with torus embeddings}, Tokyo J. Math. 5, 
37-48 (1982)

\end{thebibliography}
\end{document}